\documentclass[12pt,a4paper,twoside,final,notitlepage, leqno]{article}
\usepackage[english]{babel}
\usepackage[T1]{fontenc}
\usepackage[utf8]{inputenc}
\usepackage{epsfig, graphicx, amssymb}
\usepackage{amsmath,amsthm,epsfig,amsfonts}
\usepackage{float}
\usepackage{color}
\def\br {\break}
\def \smb {{\scriptstyle \bullet }}

\newcommand{\moneq}{\vspace*{-7pt} \begin{equation} \displaystyle }
\newcommand{\moneqstar}{\vspace*{-6pt} \begin{equation*} \displaystyle }
\newcommand{\monendstar}{\vspace*{-6pt} \end{equation*}   }
\newcommand{\monend}{\vspace*{-7pt} \end{equation}   }
\newcommand{\moneqarraystar}{ \begin{eqnarray*} \displaystyle }
\newcommand{\monendarraystar}{ \end{eqnarray*}   }
\newcommand{\R}{\mathbb{R}}

\definecolor{vertfonce}{rgb}{0.0, 0.5, 0.0}
\def\section*#1{}
\usepackage{fancyhdr}
\renewcommand{\headrulewidth}{0pt}
\begin{document}

\fancypagestyle{plain}{ \fancyfoot{} \renewcommand{\footrulewidth}{0pt}}
\fancypagestyle{plain}{ \fancyhead{} \renewcommand{\headrulewidth}{0pt}}
%%   \bibliographystyle{alpha}

%%%%%%%%%%%%%%%%%%%%%%%%%%%%%%%%%%%%%%%%%%%%%%%%%%%%%%%%%%%%%%%%%%%%%%%%%%%%%%%
~

\vskip .5 cm

\centerline {\bf \LARGE Conditions aux limites fortement non lin\'eaires}

\bigskip \smallskip 

\centerline {\bf \LARGE pour les \'equations d'Euler}

%%   \bigskip
%%   \centerline {\bf \Large pour les \'equations d'Euler de la dynamique des gaz}
\bigskip

\centerline {\bf \LARGE de la dynamique des gaz}

 \bigskip  \bigskip \bigskip \bigskip

\centerline { \large   Fran\c{c}ois Dubois$^{a}$}

\smallskip  %% \bigskip

\centerline { \it  \small
  $^a$   Aerospatiale, Division Syst\`emes Strat\'egiques et Spatiaux, Les Mureaux, France.}

%%%  \centerline { \it  \small
%%%  $^a$   Laboratoire de Math\'ematiques d'Orsay, Facult\'e des Sciences d'Orsay,}

%%%  \centerline { \it  \small   Universit\'e Paris-Saclay, France.}

%%%  \centerline { \it  \small
%%%  $^b$    Conservatoire National des Arts et M\'etiers, LMSSC laboratory,  Paris, France.}

%%%   \bigskip

\bigskip  \bigskip

\centerline {{octobre 1988} 
  {\footnote {\rm  \small $\,$
Propos\'e avec le titre ``Conditions aux limites fortement non lin\'eaires''
pour le cours CEA-EDF-INRIA 
sur les   ``M\'ethodes de diff\'e\-rences finies et \'equations hy\-perbo\-liques'', 
organis\'e par  Pierre-Louis Lions \`a l'INRIA Rocquencourt du 29 novembre au 02 d\'ecembre 1988.  \'Edition  septembre 2024.}}}
%%%    {\footnote


 {\bf AMS classification}:
%%    65Q05,   %% Difference and functional equations, recurrence relations,
76M12.   %% Finite volume methods applied to problems in fluid mechanics Recent zbMATH articles in MSC 76M12
%%  82C20.   %%%  Dynamic lattice systems (kinetic Ising, etc.) and systems on graphs

\bigskip  \bigskip  \bigskip
\noindent {\bf \large Abstract}

\noindent 
We study various formulations of the boundary conditions for the Euler equations
of gas dynamics from a mathematical and numerical point of view. In the case of one
space dimension, we recall the classical results, based on an analysis of the linearized problem.
Then we present a more recent formulation of the problem, which allows for nonlinear
effects at the boundary of the study domain. This formulation
fits naturally into a finite volume discretization, and
we present a significant one-dimensional test case.

\bigskip  \bigskip
\noindent {\bf \large R\'esum\'e}

\noindent 
Nous \'etudions diverses formulations 
des conditions aux limites pour les \'equa\-tions d'Euler
de la dynamique des gaz d'un point de vue math\'ematique et num\'erique. Dans le cas
d'une dimension d'espace, nous rappelons les r\'esultats classiques,
fond\'es sur une analyse du probl\`e\-me lin\'earis\'e et nous pr\'esentons une
formulation plus r\'ecente du probl\`eme, qui autorise la prise en compte d'effets
non lin\'eaires importants au bord du domaine d'\'etude. Cette formulation s'ins\`ere
naturellement dans une discr\'etisation par la m\'ethode des volumes finis et nous
pr\'esentons un cas test monodimensionnel significatif.

\noindent

\newpage

%%%%%%%%%%%%%%%%%%%%%%%%%%%%%%%%%%%%%%%%%%%%%%%%%%%%%%%%%%%%%%%%%%%%%%%%%%%%%%%  section 1
\noindent {\bf \large    1) \quad  Introduction} 
%%%%%%%%%%%%%%%%%%%%%%%%%%%%%%%%%%%%%%%%%%%%%%%%%%%%%%%%%%%%%%%%%%%%%%%%%%%%%%%%%%%%%%%%%%

%%%%%%%%%%%%%%%%%%%%%%%%%%%%%%%%%%       10 septembre 2024      %%%%%%%%%%%%%%%%%%%%%%%%%
\fancyhead[EC]{\sc{Fran\c{c}ois Dubois}}
\fancyhead[OC]{\sc{Conditions aux limites fortement non lin\'eaires}}
%%%%%%%%%%%%%%%%%%%%%%%%%%%%%%%%%%       10 septembre 2024      %%%%%%%%%%%%%%%%%%%%%%%%%
%%%%%%%%%%%%%%%%%%%%%%%%%%%%%%%%%%%%%%%%%%%%%  jolie numerotation des pages
\fancyfoot[C]{\oldstylenums{\thepage}}
%%%%%%%%%%%%%%%%%%%%%%%%%%%%%%%%%%%%%%%%%%%%%  fin jolie numerotation des pages

\noindent   $\bullet \qquad  \,\,\, $  
Les \'equations d'Euler de la dynamique des gaz constituent un syst\`eme hyperbolique
non lin\'eaire de lois de conservation. Nous pouvons les \'ecrire sous cette forme
dans le cas de deux dimensions spatiales par exemple~:

\smallskip \noindent  (1.1) $\qquad \displaystyle 
{{\partial W}\over{\partial t}} \,\,+\,\, {{\partial}\over{\partial x}} f(W) \,+\,
{{\partial}\over{\partial y}} g(W) \,\,= \,\, 0 \,$

\smallskip \noindent
avec des variables conservatives $\, W \,$ et des flux $\, f(W) ,\,$ $\, g(W) \,$
donn\'es par~:

\smallskip \noindent  (1.2) $\qquad \displaystyle 
W \quad \,\,\,\, = \,\, \bigl( \, \rho \,,\, \rho \,u \,,\, \rho \, v \,,\, \rho \, E
\, \bigr)^{\displaystyle \rm t} \,$ 

\smallskip \noindent  (1.3) $\qquad \displaystyle 
f(W) \,\,= \,\, \bigl( \, \rho\,u  \,,\, \rho \,u^2 \, + \,  p  \,,\, \rho \, u \, \, v
\,,\, \rho \, u \, E \,+\, p \, u  \, \bigr)^{\displaystyle \rm t} \,$ 

\smallskip \noindent  (1.4) $\qquad \displaystyle 
g(W) \,\,= \,\, \bigl( \, \rho\,v  \,,\,\rho \, u \, \, v  \,,\,  \rho \,v^2 \, + \, 
p  \,,\, \rho \, v \, E \,+\, p \, v  \, \bigr)^{\displaystyle \rm t} \,$

\smallskip \noindent
o\`u $\, \rho \,$ (respectivement $\, u ,\, v ,\, E ,\,p  $) d\'esigne la densit\'e
(repectivement les deux composantes de la vitesse, l'\'energie totale sp\'ecifique et
la pression). Il suffit de se donner la pression comme fonction des variables
conservatives pour d\'efinir compl\`etement le syst\`eme hyperbolique (1.1)~; nous
choisissons simplement une loi d'\'etat de gaz parfait polytropique~:

\smallskip \noindent  (1.5) $\qquad \displaystyle 
p \,\,= \,\, (\gamma \!-\!1) \, \rho \, \bigl( E - {1\over2} (u^2+v^2) \bigr) \,. \,$ 

\bigskip \noindent   $\bullet \qquad  \,\,\, $  
Le probl\`eme de Cauchy pos\'e sur $\R^2 \,$ pour le syst\`eme (1.1) associ\'e \`a la
condition initiale 

\smallskip \noindent  (1.6) $\qquad \displaystyle 
W(0 ,\, x ,\, y) \,\,= \,\, W_0(x,\,y) \,, \qquad (x,\,y) \in \R^2 \,$ 

\smallskip \noindent 
n'a fait, \`a notre connaissance, l'objet d'aucun r\'esultat d'existence globale en
temps, m\^eme pour des donn\'ees $\, W_0 \,$ r\'eguli\`eres. Il en est de m\^eme pour
le probl\`eme de Cauchy (1.1)(1.6) pos\'e \`a une seule dimension d'espace. Nous devons
donc aborder le probl\`eme aux limites pos\'e sur un domaine $\, \Omega \,$ de 
$ \R^2 \,$:

\setbox11=\hbox {$\displaystyle \,\, {{\partial W}\over{\partial t}} \,\,+\,\, 
{{\partial}\over{\partial x}} f(W) \,+\,
{{\partial}\over{\partial y}} g(W) \,\,= \,\, 0 \, $}
\setbox21=\hbox {$ \qquad \qquad t \geq 0 \,,\qquad (x,\,y) \in \Omega \,$}
\setbox12=\hbox {$\displaystyle \,\,W(0 ,\, x ,\, y) \,\,= \,\, W_0(x,\,y) $}
\setbox22=\hbox {$ \qquad \qquad t = 0 \,,\qquad (x,\,y) \in \Omega \,$}
\setbox13=\hbox {$\displaystyle \,\, {\rm Condition } \, \, {\rm limite} \, \bigl(
W(t ,\, x ,\, y)\bigr)  \,\,= \,\,0 $}
\setbox23=\hbox {$ \qquad \qquad t \geq 0 \,,\qquad (x,\,y) \in \partial \Omega \,$}
\setbox40= \vbox {\halign{#&# \cr \box11 & \box21\cr \box12 & \box22 \cr \box13 &
\box23 \cr}}
\setbox41= \hbox{ $\vcenter {\box40} $}
\setbox44=\hbox{\noindent  (1.7) $\quad \displaystyle \left\{ \box41 \right. $}  
\smallskip \noindent $ \box44 $

\smallskip \noindent 
avec prudence. Pourtant, c'est bien une approximation de ``la'' solution du probl\`eme
aux limites (1.7) que cherche l'ing\'enieur, dans un domaine $\, \Omega \,$ qui est
souvent {\bf non born\'e}. 

\bigskip \noindent   $\bullet \qquad  \,\,\, $  
La difficult\'e math\'ematique de l'\'etude du probl\`eme (1.1)(1.6) est li\'ee \`a
la pr\'esence d'{\bf ondes non lin\'eaires} qui imposent une \'etude en termes de
solutions faibles (voir par exemple Lax [La73] ou Smoller [Sm83]). De plus,
l'unicit\'e des solutions faibles est en g\'en\'eral en d\'efaut et une
{\bf in\'egalit\'e d'entropie} doit \^etre ajout\'ee afin d'exclure des solutions
physiquement non admissibles telles que les chocs de d\'etente par exemple
(Germain-Bader [GB53], Oleinik [Ol57], Godunov [Go61], Lax [La71]). Le probl\`emes aux
limites (1.7) n'est quant \`a lui parfaitement compris math\'ematiquement que dans le
cas {\bf lin\'eaire} (Kreiss [Kr70], Higdon [Hi86] et les r\'ef\'erences cit\'ees). 

\bigskip \noindent   $\bullet \qquad  \,\,\, $  
L'approche num\'erique pose par ailleurs un probl\`eme pratique important~: les
sch\'emas aux diff\'erences les plus classiques (Lax-Wendroff [LW60], Mac Cormack
[Mc69]) proposent un calcul {\bf centr\'e} des approximations des d\'eriv\'ees en
espace qui permettent d'incr\'ementer en temps les valeurs $\, W_{i,\,j} \,$ gr\^ace
aux valeurs $\, \, W_{i\!+\!1,\,j\!+\!1} ,\, W_{i,\,j\!+\!1} ,\,
W_{i\!-\!1,\,j\!+\!1} ,\,\dots \,\, $  situ\'ees dans un {\bf voisinage discret}
du point de grille $\, (i,\,j) .\,$ Il est donc n\'ecessaire d'introduire un ``sch\'ema
\`a la limite'' pour incr\'ementer les valeurs $\, W_{i,\,j} \,$ situ\'ees au bord du
domaine de calcul (voir par exemple Richtmyer-Morton [RM67]). La situation est m\^eme
paradoxale~: dans certaines situations physiquement bien d\'etermin\'ees (entr\'ee ou
sortie subsonique par exemple), les \'etudes lin\'eaires montrent que le probl\`eme
est bien pos\'e avec moins de conditions limites que le syst\`eme (1.1) ne compte
d'\'equations (voir par exemple Oliger-Sundstr\"om [OS78], Yee-Beam-Warming [YBW82] ou
Gustafsson [Gu85]). Les conditions aux limites suppl\'ementaires doivent \'egalement
\^etre telles que le sch\'ema {\bf global} ({\it i.e.} le sch\'ema \`a l'int\'erieur et le
sch\'ema \`a la limite) reste stable. Le travail classique de Gustafsson, Kreiss et
Sundstr\"om [GKS72] a permis de d\'evelopper une m\'ethode d'analyse de la stabilit\'e
des sch\'emas aux diff\'erences finies pour des syst\`emes hyperboliques
{\bf lin\'eaires}. 

\bigskip \noindent   $\bullet \qquad  \,\,\, $  
Dans ces notes de cours, nous nous limiterons essentiellement \`a des pro\-bl\`emes
{\bf monodimensionnels}, ce qui correspond pour le probl\`eme aux limi\-tes (1.7) \`a une
analyse dans la direction normale \`a la fronti\`ere. Dans une premi\`ere partie,
nous \'etudions le probl\`eme continu. Au cours d'une seconde partie, nous nous
int\'eressons aux sch\'emas num\'eriques existants pour l'\'ecriture de conditions
aux limites lin\'eaires ou non. Enfin nous proposons une m\'ethode de type volumes
finis qui \'etend tr\`es simplement l'approche originale de Godunov [Go59] (voir
aussi Godunov {\it et al.}  [GZIKP79]) et permet la prise en compte de fortes
non-lin\'earit\'es \`a la fronti\`ere du domaine de calcul.

\bigskip \bigskip 
%%%%%%%%%%%%%%%%%%%%%%%%%%%%%%%%%%%%%%%%%%%%%%%%%%%%%%%%%%%%%%%%%%%%%%%%%%%%%%%  section 2
\noindent {\bf \large    2) \quad  Etude du probl\`eme continu} 
%%%%%%%%%%%%%%%%%%%%%%%%%%%%%%%%%%%%%%%%%%%%%%%%%%%%%%%%%%%%%%%%%%%%%%%%%%%%%%%%%%%%%%%%%%

\smallskip \noindent {\bf 2.1)  \quad  	 Quelques rappels \'el\'ementaires}

\smallskip 
\noindent   $\bullet \qquad  \,\,\, $  
Nous r\'e\'ecrivons le syst\`eme (1.1) des \'equations d'Euler dans le cas d'une
seule dimension d'espace. La composante $\, v \,$ de la vitesse est  identiquement
nulle et l'on a simplement~:

\smallskip \noindent  (2.1) $\qquad \displaystyle 
{{\partial W}\over{\partial t}} \,\,+\,\, {{\partial}\over{\partial x}} f(W)  \,\,=
\,\, 0 \,$

\smallskip \noindent
avec 

\smallskip \noindent  (2.2) $\qquad \displaystyle 
W \quad \,\,\,\, = \,\, \bigl( \, \rho \,,\, \rho \,u \,,\, \rho \, E \,
\bigr)^{\displaystyle \rm t} \,$

%% \smallskip
\noindent  (2.3) $\qquad \displaystyle 
f(W) \,\,= \,\, \bigl( \, \rho\,u  \,,\, \rho \,u^2 \, + \,  p  \,,\,  \rho \, u \, E
\,+\, p \, u  \, \bigr)^{\displaystyle \rm t} \,. \,$ 

\smallskip \noindent
La pression est calcul\'ee gr\^ace \`a la loi d'\'etat (1.5). Il est utile d'\'ecrire
(2.1) sous  forme {\bf non conservative} et nous choisissons pour cela les variables
$\, V \,$ d\'efinies par~: 

\smallskip \noindent  (2.4) $\qquad \displaystyle 
V \,\,= \,\, \bigl( \, \rho \,,\, u \,,\, S \,  \bigr)^{\displaystyle \rm t} \,. \,$ 

\smallskip \noindent 
L'entropie sp\'ecifique adimensionnalis\'ee $\, S \,$ est reli\'ee aux autres
variables thermodynamiques $\, \rho ,\, p \,$ par la relation~: 

\smallskip \noindent  (2.5) $\qquad \displaystyle 
p \,\,= \,\, S \, \rho^{\gamma} \,$

\smallskip \noindent
(voir les textes classiques tels que Courant-Friedrichs [CF48] ou Landau-Lifchitz
[LL54]) et le syst\`eme (2.1) s'\'ecrit sous forme \'equivalente 

\smallskip \noindent  (2.6) $\qquad \displaystyle 
{{\partial V}\over{\partial t}} \,\,+\,\, A(V) \, {{\partial V}\over{\partial x}}
\,\,= \,\, 0 \,$

\smallskip \noindent
lorsque la solution $\, V(t,\,x) \,$ est r\'eguli\`ere. Nous avons~:

\smallskip \noindent  (2.7) $\qquad \displaystyle 
A(V) \,= \, \begin{pmatrix} u & \rho & 0 \\   {{1}\over{\rho}} \,
{{\partial p}\over{\partial \rho}}  & u &  {{1}\over{\rho}} \,
{{\partial p}\over{\partial S}} \\ 0 & 0 & u \end{pmatrix} $. 

\smallskip \noindent 
La c\'el\'erit\'e du son $\, \, c ,\,$ d\'efinie par

\smallskip \noindent  (2.8) $\qquad \displaystyle 
c \,\,= \,\, \sqrt{ {{\partial p}\over{\partial \rho}}(\rho,\,S) \,} \,\,= \,\,
\sqrt{{{\gamma \, p }\over{\rho}} \, } \,$ 

\smallskip \noindent 
intervient dans le calcul des valeurs propres $\,\, \lambda_{j}(W) \,\,$ de la
matrics $\, A(V) \,$: 

\smallskip \noindent  (2.9) $\qquad \displaystyle 
\lambda_{1}(W) \, \equiv \, u-c  \quad  < \quad  \lambda_{2}(W) \,\equiv \, u
\quad  < \quad  \lambda_{3}(W) \, \equiv \, u+c \,.\,  $

\smallskip \noindent 
Les vecteurs propres associ\'es \`a ces valeurs propres se calculent facilement~:

\smallskip \noindent  (2.10) $\qquad \displaystyle 
r_{1}(V) \,=\, \left( \,\, \begin{array}{c} 
  \rho \\  \vspace{.1 cm}   \\ -c \\ \vspace{.1 cm}   \\ 0 \end{array} \right) \,\,; \quad 
r_{2}(V) \,=\,\left( \,\, \begin{array}{c} 
  {\partial p} /{\partial S} \\  \vspace{.1 cm}   \\ 0 \\ \vspace{.1 cm}   \\  -c^2  \end{array} \right) \,\,; \quad 
 r_{3}(V) \,=\, \left( \,\, \begin{array}{c} 
  \rho \\  \vspace{.1 cm}   \\ c \\ \vspace{.1 cm}   \\ 0 \end{array} \right) $.

\bigskip \noindent   $\bullet \qquad  \,\,\, $  
Le {\bf probl\`eme de Riemann} $\,\, R(W_g,\,W_d) \,\,$ associ\'e au syst\`eme (2.1)
est un probl\`eme de Cauchy particulier~; la condition initiale est compos\'ee de
deux \'etats constants~:

\setbox11=\hbox {$\displaystyle \,\,W_g \qquad \qquad x < 0  \, $}
\setbox12=\hbox {$\displaystyle \,\,W_d \qquad \qquad x > 0  \,. \, $}
\setbox40= \vbox {\halign{#&# \cr \box11 \cr  \box12 \cr  }}
\setbox41= \hbox{ $\vcenter {\box40} $}
\setbox44=\hbox{\noindent  (2.11) $\qquad \displaystyle W_0(x) \,\,= \,\, 
\left\{ \box41 \right. $}  
\smallskip \noindent $ \box44 $

\smallskip \noindent
La solution entropique de ce probl\`eme est constitu\'ee d'ondes de choc, d'ondes de
d\'etente et d'une discontinuit\'e de contact, s\'epar\'es par au plus deux \'etats
constants (voir par exemple Courant-Friedrichs [CF48], Landau-Lifchitz [LL54],
Godunov  {\it et al.}   [GZIKP79]). Rappelons que le long d'une d\'etente, les invariants de
Riemann associ\'es sont constants, {\it i.e.} 

\smallskip \noindent  (2.12) $\qquad \displaystyle 
w_1^1 \,\,= \,\, S \,\, ; \qquad w_2^1 \,\,= \,\, u \,+\, {{2}\over{\gamma\!-\!1}} \,c
\hfill 1 $-d\'etente

\smallskip \noindent  (2.13) $\qquad \displaystyle 
w_1^2 \,\,= \,\, u \,\, ; \qquad w_2^2 \,\,= \,\, p  
\hfill  2 $ -discontinuit\'e de contact 

\smallskip \noindent  (2.14) $\qquad \displaystyle 
w_1^3 \,\,= \,\, S \,\, ; \qquad w_2^3 \,\,= \,\, u \,-\, {{2}\over{\gamma\!-\!1}} \,c
\hfill  3 $ -d\'etente.

\smallskip \noindent 
Rappelons \'egalement que le long d'une telle d\'etente, $\, \, W \,=\, W(x/t) \,\,$
est une solution autosemblable qui satisfait \`a l'\'equation diff\'erentielle~:

\smallskip \noindent  (2.15) $\qquad \displaystyle 
{{ {\rm d}} \over{ {\rm d} \, \Bigl( {{\displaystyle  x}\over{\displaystyle  t}} 
\Bigr) }} \, \biggl( W \Bigl( {{\displaystyle  x}\over{\displaystyle  t}} 
\Bigr) \biggr) \,\,\,=\,\, \, r_{j} \,  \biggl(   W\Bigl( {{\displaystyle 
x}\over{\displaystyle  t}}  \Bigr) \biggr) \, \qquad \qquad j \,= \, 1 \,\, {\rm ou }
\,\, 3 \,. \,$ 

\smallskip \noindent
A travers une onde de choc ou une discontinuit\'e de contact de c\'el\'erit\'e  $\,
\sigma , \,$ on a les relations de Rankine-Hugoniot~:

\smallskip \noindent  (2.16) $\qquad \displaystyle 
\bigl[ \, f(W) \, \bigr] \,\,= \,\, \sigma \,  \bigl[ \, W \, \bigr] \,.\, $ 

\smallskip \noindent 
Rappelons que le champ num\'ero 2 (lin\'eairement d\'eg\'en\'er\'e) peut \^etre vu
\`a la fois comme un choc et une d\'etente. La solution du probl\`eme de Riemann
consiste \`a construire trois $j$-courbes $\, \, U_{j}(W) \,\,$ dans l'espace des
phases (ou des \'etats) au voisinage de chaque \'etat $\, W \,$~:

\setbox11=\hbox {$\displaystyle \,\,W \in U_{j}(W) \,$ }
\setbox12=\hbox {$\displaystyle \,\,\forall \, W' \in U_{j}(W), \,\, $ la solution 
de $\, R(W,\,W') \,$ ou de $\, R(W',\,W) \,$ }
\setbox13=\hbox { \qquad \qquad \qquad  \qquad  est une onde simple d\'ecrite plus
haut.  }
\setbox40= \vbox {\halign{#&# \cr \box11 \cr  \box12  \cr \box13  \cr }}
\setbox41= \hbox{ $\vcenter {\box40} $}
\setbox44=\hbox{\noindent  (2.17) $\qquad \displaystyle \left\{ \box41 \right. $}  
\smallskip \noindent $ \box44 $

\smallskip \noindent
Les \'etats interm\'ediaires $\, W_1 \,$ et $\, W_2 \,$ sont alors d\'efinis de sorte
que 

\smallskip \noindent  (2.18) $\qquad \displaystyle 
W_1 \in U_{1}(W_g) \,\,; \qquad W_2 \in U_{2}(W_1) \,\,; \qquad W_d \in 
U_{3}(W_2) \,.\,$ 

\smallskip \noindent 
Pour plus de d\'etails concernant les courbes $\, \, U_{j}({\displaystyle
\bullet}),\,$ nous renvoyons \`a Lax [La73] ou Smoller [Sm83]. 

\bigskip \noindent {\bf 2.2)  \quad  	 Syst\`eme des \'equations d'Euler
lin\'earis\'ees}

\noindent   $\bullet \qquad  \,\,\, $  
La plupart des r\'esultats existants concernent les \'equations hyperboliques
lin\'eaires. Nous lin\'earisons donc les \'equations d'Euler autour d'un {\bf \'etat
constant} $\,\, {\overline W} \,\,$ en posant~:

\smallskip \noindent  (2.19) $\qquad \displaystyle 
W \,\,= \,\,  {\overline W} \,+\, W' \,$ 

\smallskip \noindent
et en n\'egligeant les termes du second ordre en $\, W' \,$ dans le syst\`eme obtenu.
Les r\'esultats les plus simples r\'esultent de la forme (2.6) des \'equations et
l'on obtient~:

\smallskip \noindent  (2.20) $\qquad \displaystyle 
{{\partial V'}\over{\partial t}} \,+\, A( {\overline V} ) \, {{\partial
V'}\over{\partial x}} \,\,= \,\, B({\overline V} ,\, V') \,.\,$ 

\smallskip \noindent 
On peut alors poursuivre l'\'etude de ce syst\`eme en diagonalisant la matrice $\, A(
{\overline V} ) .\,$ Nous notons $\, \varphi_{j} \,$ les coordonn\'ees de la
perturbation $\, V' \,$ dans la base (fixe) des vecteurs propres $\, \, r_{j}(
{\overline V} )\,$~:

\smallskip \noindent  (2.21) $\qquad \displaystyle 
V' \,\,= \,\, \sum_{j = 1}^{3} \, \varphi_{j} \,\, r_{j}({\overline V} )\,$

\smallskip \noindent 
et l'on a~: 

\setbox11=\hbox {$\displaystyle \,\,\varphi_{1} \,\, = \,\, {{1}\over{2 \,
{\overline \rho} \, {\overline c}^2 }} \, \bigl( \, p' \,-\, {\overline \rho} \,
{\overline c} \, u' \, \bigr) $ }
\setbox12=\hbox {$\displaystyle \,\,\varphi_{2} \,\, = \,\, \,\,\, 
-{{1}\over{{\overline c}^2 }} \, S' $ }
\setbox13=\hbox {$\displaystyle \,\,\varphi_{3} \,\, = \,\, {{1}\over{2 \,
{\overline \rho} \, {\overline c}^2 }} \, \bigl( \, p' \,+\, {\overline \rho} \,
{\overline c} \, u' \, \bigr) \,. \, $ }
\setbox40= \vbox {\halign{#&# \cr \box11 \cr  \box12  \cr \box13  \cr }}
\setbox41= \hbox{ $\vcenter {\box40} $}
\setbox44=\hbox{\noindent  (2.22) $\qquad \displaystyle \left\{ \box41 \right. $}  
\smallskip \noindent $ \box44 $

\bigskip \noindent   $\bullet \qquad  \,\,\, $  
La diff\'erence de pression $\,\,p' \,\,$ est donn\'ee selon la relation
lin\'earis\'ee~: 

\smallskip \noindent  (2.23) $\qquad \displaystyle 
p' \,\,\equiv \,\, p \,-   {\overline p} \,\,= \,\, {{ \partial p}\over{\partial S}}(
{\overline W}) \, S' \,+\,  {\overline c}^2 \, \rho'  \,.\,$ 

\smallskip \noindent  
Le changement de variables $\,\, V' \longmapsto \varphi \,\,$ permet de d\'ecoupler
le membre de gauche de la relation (2.20) sous la forme de trois \'equations
d'advection~:

\smallskip \noindent  (2.24) $\qquad \displaystyle 
{{\partial \varphi}\over{\partial t}} \,+\, \Lambda( {\overline W})\, {{\partial
\varphi}\over{\partial x}} \,\,= \,\, C( {\overline W},\, \varphi) \,\,$ 

\smallskip \noindent 
o\`u $\,\, \Lambda( {\overline W})\,=\ {\rm diag} ({\overline u}-{\overline c} \,,\,
{\overline u} \,,\,{\overline u}+{\overline c}). \,$ Le couplage entre les
composantes de $\, \, \varphi \,\,$ est uniquement r\'ealis\'e par le second membre
$\,\, C( {\overline W},\, \varphi) ,\,$ qui est un op\'erateur non diff\'erentiel. 
Le syst\`eme (2.24) est appel\'e dans la suite ``syst\`eme des \'equations d'Euler
sous forme caract\'eristique'' et les variables $\,\, \varphi \,\,$ d\'efinies en
(2.22) sont les ``variables caract\'eristiques''. L'int\'er\^et essentiel de cette
d\'emarche est qu'on a ramen\'e, dans le cas lin\'earis\'e, le syst\`eme des
\'equations d'Euler \`a une forme classique. 

\bigskip \noindent {\bf 2.3)  \quad  	 Probl\`eme aux limites pour un
syst\`eme hyperbolique lin\'eaire}

\noindent   $\bullet \qquad  \,\,\, $  
Dans un article devenu classique, Kreiss [Kr70] introduit une notion de ``probl\`eme
bien pos\'e'' pour l'\'etude du probl\`eme aux limites associ\'e au syst\`eme
(2.24) dans le ``quart d'espace'' $\,\, t\geq 0 \,,\, x \geq 0 .\,\,$ Nous notons $\,\,
\Lambda^{I} \,$ (respectivement $\, \Lambda^{II}  $) la matrice diagonale obtenue
\`a partir de $\,\, \Lambda({\overline W})\,\,$ en ne conservant que les valeurs
propres positives (respectivement n\'egatives), donc nous supposons ici $\,\,
{\overline u} \not= 0.\,\,$ Nous d\'ecomposons \'egalement les variables
caract\'eristiques $\,\, \varphi \,\,$ sous la forme 

\smallskip \noindent  (2.25) $\qquad \displaystyle 
\varphi \,\, = \,\, \bigl( \, \varphi^I \,,\, \varphi^{II} \, \bigr) \,$ 

\smallskip \noindent 
de fa\c{c}on \`a utiliser les notations de Kreiss. La condition limite en $\, x=0 \,$
est \'ecrite sous la forme 

\smallskip \noindent  (2.26) $\qquad \displaystyle 
\varphi^I  \,\,= \,\, \Sigma  \, \varphi^{II}  \,+\, g \qquad \qquad t \geq 0 \,,\,\,\,
x \geq 0 \,$ 

\smallskip \noindent
et la condition initiale est simplement

\smallskip \noindent  (2.27) $\qquad \displaystyle 
\varphi(0,\,x) \,\, = \,\, \varphi_0(x) \qquad  \qquad  t = 0 \,,\,\,\, x \geq 0 \,.
\,$ 

\smallskip \noindent
Le syst\`eme (2.24) associ\'e \`a la condition (2.26) et \`a la condition initiale
(2.27) (l'``IVBP'' (2.24)(2.26)(2.27) de fa\c{c}on plus concise) est alors {\bf bien 
pos\'e} dans $\, {\rm L}^2 ,\,$ au sens de Kreiss. 

\bigskip \noindent   $\bullet \qquad  \,\,\, $  
La condition limite (2.26) s'interpr\`ete en termes de {\bf directions
carac\-t\'eristiques}~: le champ $\, \varphi \,$ le long des caract\'eristiques
entrantes est une  {\bf fonction affine} des composantes le long des
caract\'eristiques sortant du domaine d'\'etude.  Nous insistons sur le fait que le
r\'esultat pr\'ec\'edent peut s'\'etendre dans diverses directions. En particulier le
cas multidimensionnel (beaucoup plus complexe) peut \^etre abord\'e \`a l'aide des
m\^emes concepts~: la direction $\,x\,$ doit alors \^etre rempla\c{c}\'ee par la
normale au domaine $\, \Omega \,$ (voir par exemple Higdon [Hi86]). Dans le cas o\`u
la vitesse de r\'ef\'erence $\,\,  {\overline u}\,\,$ est nulle (fronti\`ere
caract\'eristique) les r\'esultats de Kreiss ont \'et\'e \'etendus par Majda-Osher
[MO75]. 

\bigskip \noindent   $\bullet \qquad  \,\,\, $  
Pour la mise en \oe uvre pratique de la condition limite (2.26), on distingue
habituellement quatre cas selon que la vitesse $\,\,  {\overline u}\,\,$ est
positive (entr\'ee) ou n\'egative (sortie), de module sup\'erieur \`a la
c\'el\'erit\'e du son (supersonique) ou inf\'erieur (subsonique). 

(i) \qquad Entr\'ee supersonique $\,\,( {\overline u} >  {\overline c}). \,$ 

\noindent La composante ``sortante'' de $\, \varphi ,\,$ \`a savoir $\,\, \varphi^{II},
\,\,$ est nulle et la relation (2.26) revient \`a se donner toutes les composantes de
$\, \varphi .\,$

(ii) \qquad Entr\'ee subsonique $\,\,( 0 < {\overline u} <  {\overline c}). \,$ 

\noindent On dispose de deux caract\'eristiques entrantes $\,\, \varphi^I \,\,$ et
d'une caract\'eristique sortante $\,\, \varphi^{II} .\,$ Le probl\`eme lin\'earis\'e
est donc bien pos\'e lorsqu'on se donne l'un  des couples suivants 
({\it c.f.}  Oliger-Sundstr\"om [OS78] ou Yee-Beam-Warming 
  [YBW82])~: (densit\'e, pression),
(vitesse, pression) ou (enthalpie, entropie). 

(iii) \qquad Sortie subsonique $\,\,(- {\overline c} < {\overline u} <  0) .\,$ 

\noindent  Une seule caract\'eristique entre dans le domaine de calcul et deux sont
sortantes. Il est classique de se donner la presion ou la vitesse de sortie.
Remarquons que le choix d'une pression impos\'ee $\,\,p = {\overline p}\,\,$
s'\'ecrit apr\`es lin\'earisation autour de $\,\, {\overline V} \,=\, ( {\overline
\rho} ,\,  {\overline u}  ,\,  {\overline p}  )^{\displaystyle \rm t} \,$:

\smallskip \noindent  (2.28) $\qquad \displaystyle 
p' \,\,= \,\, 0 \,$

\smallskip \noindent 
ce qui revient, dans la relation (2.26), \`a prendre 

\smallskip \noindent  (2.29) $\qquad \displaystyle 
\Sigma \,\,= \,\, \bigl( \, -1 \,,\, 0 \, \bigr) \,,\, \qquad g \,\,= \,\, 0 \,.\,$ 

\newpage
(iv) \qquad Sortie supersonique $\,\,(  {\overline u} <   -{\overline c}) .\,$ 

\noindent  Toutes les caract\'eristiques sont sortantes $\,\, (\varphi^{I} = 0) ,\,$
donc aucune information n'est contenue dans la relation (2.26)~; aucune ``condition
analytique'' n'est n\'ecessaire dans ce cas. 

\bigskip \noindent   $\bullet \qquad  \,\,\, $  
Aux quatre cas pr\'ec\'edents, il convient de rajouter le cas singulier o\`u
$\,\, {\overline u} = 0 \,$ qui correspond physiquement \`a une paroi solide. On montre
(Oliger-Sundstr\"om [OS78] par exemple) que le probl\`eme aux limites
(2.24)(2.27)(2.30) est bien pos\'e dans $\, {\rm L}^2 \,$ avec la condition naturelle 

\smallskip \noindent  (2.30) $\qquad \displaystyle 
u' \,\,= \,\, 0 \,$ 

\smallskip \noindent 
qui revient \`a imposer l'imperm\'eabilit\'e de la paroi pour les \'equations
lin\'earis\'ees. Remarquons que (2.30) peut \'egalement s'\'ecrire sous la forme
(2.26) avec le choix 

\smallskip \noindent  (2.31) $\qquad \displaystyle
\Sigma \,\,= \,\, \bigl( \, 1 \,,\, 0 \, \bigr) \,,\, \qquad g \,\,= \,\, 0 \,.\,$ 

\bigskip \noindent {\bf 2.4)  \quad  	Probl\`eme aux limites dans le cas non lin\'eaire}

\noindent   $\bullet \qquad  \,\,\, $  
Les \'equations d'Euler sont abord\'ees th\'eoriquement dans le cas monodimensionnel
\`a l'aide de deux approches math\'ematiques~: la m\'ethode de Glimm (Glimm [Gl65])
 et la compacit\'e
par compensation (DiPerna [DP83]). La compacit\'e par compensation  permet  \`a
notre connaissance de n'aborder que des syst\`emes hyperboliques de deux \'equations
seulement et aucune \'etude d'un pro\-bl\`eme aux limites fond\'ee sur cette technique
n'a encore \'et\'e propos\'ee. La m\'ethode de Glimm permet de prouver l'existence de
solutions entropiques pour le pro\-bl\`eme de Cauchy (1.1)(1.6)  pos\'e sur $\, \R
\,$({\it i.e.} $\, x \in \R ,\, t \geq 0 $) pour un syst\`eme hyperbolique quelconque
lorsque la condition initiale est proche d'un \'etat $\,  {\overline W} \,$ fix\'e.
Il est donc naturel de faire la m\^eme hypoth\`ese lorsqu'on s'int\'eresse au
pro\-bl\`eme aux limites. L'analyse du nombre de conditions scalaires \`a imposer,
d\'evelopp\'ee au paragraphe pr\'ec\'edent, a permis de formuler de fa\c{c}on
raisonnable les conditions aux limites. Cette approche a \'et\'e suivie par
Nishida-Smoller [NS77] et Liu [Li77] lors de leur \'etude du $p$-syst\`eme de la
dynamique des gaz isentropiques. Ces auteurs montrent que, associ\'ee \`a une
condition limite sur la pression ou la vitesse, la m\'ethode de Glimm converge lorque
le pas du  maillage tend vers z\'ero. Le choix d'une condition limite, \'ecrite sous
forme {\bf forte}, r\'esulte de l'analyse du syst\`eme lin\'earis\'e. On peut donc la
qualifier de ``faiblement non lin\'eaire''. 

\bigskip \noindent   $\bullet \qquad  \,\,\, $  
Le cas d'une inconnue $\, w \,$ {\bf scalaire} est tr\`es int\'eressant puisqu'on
dispose d'un th\'eor\`eme d'existence et d'unicit\'e pour le probl\`eme de Cauchy
(pos\'e dans $\, \R \,$ ou $\, \R^2 $) (Kru$\breve {\rm z} $kov [Kv70]). Le
probl\`eme aux limites correspondant a \'et\'e abord\'e par Bardos-Leroux-N\'ed\'elec
[BLN79]. Cette \'etude a \'et\'e g\'en\'eralis\'ee aux syst\`emes, mais en se
restreignant au cas monodimensionnel dans notre travail avec Philippe Le Floch [DL88]. 
De fa\c{c}on pr\'ecise, nous montrons que sous des hypoth\`eses raisonnables de convergence,
la limite $\,  w \,$ de l'approximation visqueuse $\, w^{\epsilon} \,$ solution du
probl\`eme parabolique suivant~: 

\setbox11=\hbox{$\displaystyle \,\,{{ \partial w^{\epsilon}}\over{\partial t}}
\,+\,{{ \partial }\over{\partial x}}  \, f\bigl( w^{\epsilon} \bigr) \,\,=\,\,
\epsilon \,  {{ \partial^2  w^{\epsilon} }\over{\partial x^2}}  \,$  }
\setbox21=\hbox {$ \qquad \qquad x > 0 \,,\qquad t > 0  \,$}
\setbox12=\hbox {$\displaystyle \,\,w^{\epsilon}(0,\,x) \,\,= \,\, v_0(x) \, $}
\setbox22=\hbox {$ \qquad \qquad x > 0   \,$}
\setbox13=\hbox {$\displaystyle \,\,w^{\epsilon}(t,\,0) \,\,= \,\,u_0(t)  $}
\setbox23=\hbox {$ \qquad \qquad t > 0   \,$}
\setbox40= \vbox {\halign{#&# \cr \box11 & \box21\cr \box12 & \box22 \cr \box13 &
\box23 \cr}}
\setbox41= \hbox{ $\vcenter {\box40} $}
\setbox44=\hbox{\noindent  (2.32) $\qquad \displaystyle \left\{ \box41 \right. $}  
\smallskip \noindent $ \box44 $

\smallskip \noindent
v\'erifie au bord du domaine une {\bf in\'egalit\'e d'entropie}, pour tout couple $\,
(\eta,\, \xi) \,$ d'entropie-flux au sens de Lax [La71]~:

\smallskip \noindent  (2.33) $ \quad \displaystyle
\xi \bigl( w(t,\,0^+) \bigr)  - \xi \bigl( w_0(t)  \bigr) \,-\, {\rm d}\eta
\bigl( w_0(t)  \bigr) \,{\scriptstyle \bullet}\, \Bigl(  f \bigl(  w(t,\,0^+) \bigr)
 -   f \bigl(  w_0(t) \bigr)\Bigr) \, \leq \, 0 \,$ 

\smallskip \noindent  
avec $\,\,  w(t,\,0^+) \,=\, \displaystyle \lim_{x \rightarrow  0} \, 
w(t,\,x) \,.$ 

\bigskip \noindent   $\bullet \qquad  \,\,\, $  
Cette in\'egalit\'e g\'en\'eralise au cas des syst\`emes hyperboliques l'approche
propos\'ee initialement par Bardos-Leroux-N\'ed\'elec [BLN79] pour les lois de
conservation scalaires et dans le cas o\`u $\, \eta \,$ est une entropie de 
Kru$\breve {\rm z} $kov~: 
$\,\, \eta(w) \,=\, $ $ \mid w-k \mid . \, $ Elle a \'egalement \'et\'e
obtenue ind\'ependemment dans le cas g\'en\'eral par Audounet [Au84] et Mazet  {\it et al.} 
[MBGB87] \`a partir d'une formulation variationnelle entropique des \'equations
d'Euler. Dans le cas d'une \'equation scalaire avec poids, nous renvoyons \`a Le
Floch et N\'ed\'elec [LN88]. Nous proposons ici, comme dans  [DL88] de 
{\bf d\'efinir} la condition limite par la relation 

\smallskip \noindent  (2.34) $\qquad \displaystyle
w(t,\,0^+) \in {\cal E} \bigl( w_0(t) \bigr) \,\,, \qquad \qquad t > 0 \,$ 

\smallskip \noindent
o\`u l'ensemble limite $\, {\cal E} \bigl( w_0 \bigr) \,$ associ\'e \`a la
``condition limite'' $\, W_0 \,$ est d\'efini par~:

\setbox11=\hbox{$\displaystyle w \,  / \,\,    \xi \bigl( w \bigr)- \xi \bigl( w_0
\bigr) - {\rm d}\eta \bigl( w_0  \bigr) \,{\scriptstyle \bullet}\, \Bigl(  f
\bigl( w \bigr) -  f \bigl(  w_0\bigr) \Bigr) \, \leq \, 0 \,$  }
\setbox12=\hbox {$\displaystyle \qquad  \forall \, (\eta ,\, \xi ) \,\,$ couple
entropie-flux au sens de Lax  }
\setbox40= \vbox {\halign{#&# \cr \box11 \cr \box12 \cr}}
\setbox41= \hbox{ $\vcenter {\box40} $}
\setbox44=\hbox{\noindent  (2.35) $\quad \displaystyle  {\cal E} \bigl( w_0 \bigr) 
\,\,= \,\,  \left\{ \box41 \right\} \, ; $}  
\smallskip \noindent $ \box44 $

\smallskip \noindent
cette notion de condition \`a la limite constitue une extension de la condition de
Dirichlet habituelle. 

\bigskip \noindent   $\bullet \qquad  \,\,\, $  
Pour un syst\`eme hyperbolique lin\'eaire, le choix d'une entropie particuli\`ere
permet de montrer que l'in\'egalit\'e d'entropie \`a la limite (2.33) est
\'equivalente \`a la condition classique (2.26) avec le choix $\,\, \Sigma = 0 \,\,$
(pas de r\'eflexion d'onde). Nous le d\'etaillons pour l'\'equation d'{\bf advection}

\smallskip \noindent  (2.36) $\qquad \displaystyle
{{\partial w}\over{\partial t}} \,+\, a \, {{\partial w}\over{\partial x}} \,\,= \,\,
0 \,, \qquad \qquad w(t,\,x) \in \R \,$ 

\smallskip \noindent
qui admet le couple entropie-flux $\,\, (\eta,\, \xi) \,=\, (w^2 ,\, a \, w^2). \,$
L'in\'egalit\'e (2.33) s'\'ecrit alors dans ce cas particulier 

\smallskip \noindent  (2.37) $\qquad \displaystyle
a \, \bigl( w - w_0 \bigr)^2 \,\, \leq \,\, 0 \,$ 

\smallskip \noindent
et nous pouvons en d\'eduire facilement la discussion classique sur le nombre de
conditions aux limites du paragraphe 2.

\bigskip \noindent   $\bullet \qquad  \,\,\, $  
Le cas de l'\'equation de Burgers $\,\, \bigl( f(w) \,=\, w^2 / 2 \,,\,\, w \in \R
\bigr) \,\,$ permet de calculer compl\`etement l'ensemble  $\,\, {\cal E} \bigl( w_0
\bigr) \,\,$ (voir aussi Le Floch [LF88]) et l'on a~: 

\setbox11=\hbox{$\displaystyle \,\, ] -\infty \,,\, -w_0  ] \, \cup \, \{w_0 \} \qquad
\qquad w_0 \geq 0 \,  $  }
\setbox12=\hbox {$\displaystyle \,\, ] -\infty \,,\, 0  ] \qquad  \qquad \quad
\qquad \qquad w_0 \leq 0 \,. \,   $  }
\setbox40= \vbox {\halign{#&# \cr \box11 \cr \box12 \cr}}
\setbox41= \hbox{ $\vcenter {\box40} $}
\setbox44=\hbox{\noindent  (2.38) $\qquad \displaystyle  {\cal E} \bigl( w_0 \bigr) 
\,\,= \,\,  \left\{ \box41 \right. \,  $}  
\smallskip \noindent $ \box44 $

\bigskip  \noindent 
Nous lisons la relation (2.38) de la fa\c{c}on suivante~: lorsque $\, w_0 \,$ est
strictement positif (entr\'ee ``supersonique'' en m\'ecanique de Burgers) et l'\'etat
limite $\, w \,$ ``proche'' de celui-ci, alors $\, w \,$ est \'egal \`a $\, w_0 \,$ et
l'\'etude lin\'earis\'ee du probl\`eme \`a la limite (qui revient dans ce cas \`a
\'etudier l'\'equation d'advection) fournit la bonne condition. Lorsque $\, w \,$ est
``assez loin'' ({\it i.e.}\br
$ \mid w - w_0 \mid \, > \, 2 \, w_0  ) \,\,$ on peut accepter
\`a la limite des \'etats ``sortants'' et la condition d'entropie (2.34)(2.35) 
s'av\`ere ``fortement non lin\'eaire''.

\newpage

%%%%%%%%%%%%%%%%%%%%%%%%%%%%%%%%      figure 2-1    %%%%%%%%%%%%%%%%%%%%%%%%%%%%%%%% 
%%%%%  \bigskip 
%%%%%   \centerline {  \epsfysize=9,5cm  \epsfbox  {../cnslim/Cnslim.fig1.epsf} } 
%%%%%  \smallskip \smallskip
%%%%%  \centerline { {\bf Figure 1}	\quad 	 Ensemble limite  $\,\, {\cal V} \bigl( W_0
%%%%%  \bigr)  \,\,$  [repr\'esent\'e en trait fort }
%%%%%  \centerline { et \`a l'aide des  zones gris\'ees] pour les \'equations d'Euler-Saint
%%%%%  Venant }
%%%%%  \centerline { de la dynamique des gaz dans l'approximation isentropique. }
%%%%%  \centerline {  Cas d'une entr\'ee supersonique. } 
%%%%%%%%%%%%%%%%%%%%%%%%%%%%%%%%%%%%%%%%%%%%%%%%%%%%%%%%%%%%%%%%%%%%%%%%%%%%%%%%%%%
~

\vskip -1.2 cm 
%%%%%%%%%%%%%%%%%%%%%%%%%%%%%%%%%%%%%%%%%%%%%%%%%%%%%%%%%%%%%%%%%%%%%%%%%%%%%%%%%%% figure  1
\centerline  {\includegraphics[width=.55\textwidth]   {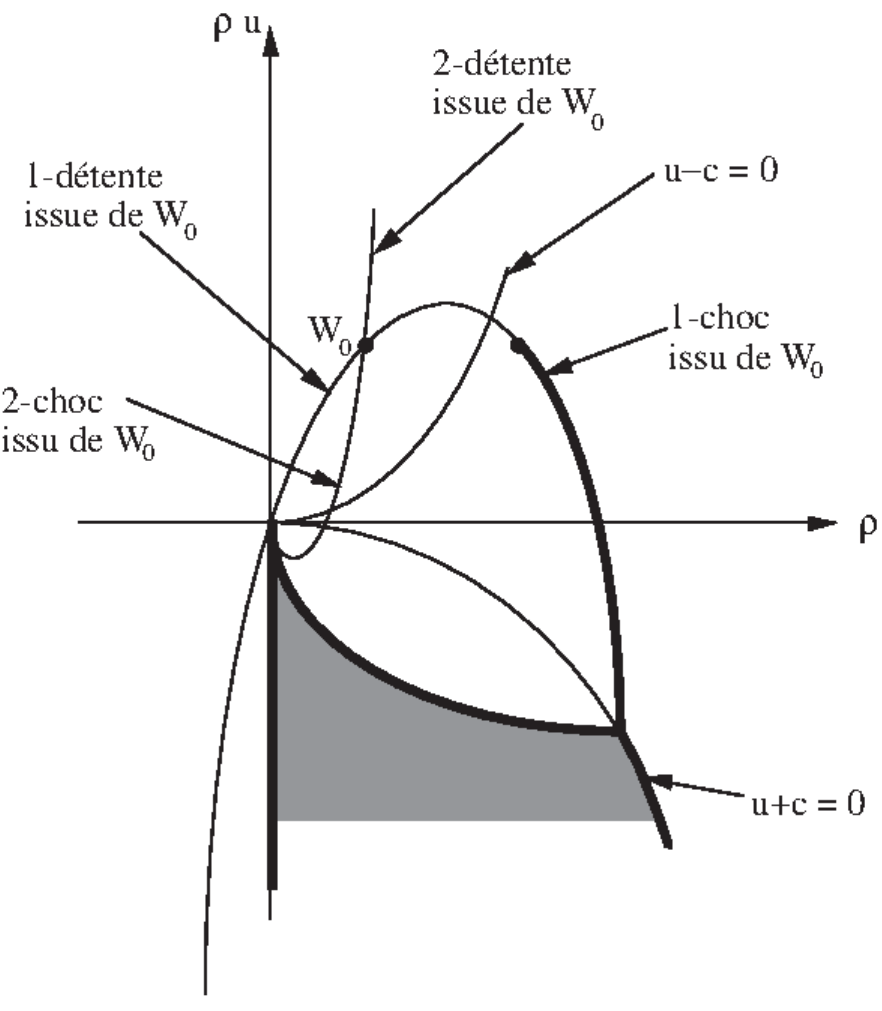}}
\smallskip  
\noindent {{\bf Figure 1}. Ensemble limite  $\, {\cal V} \bigl( W_0 \bigr)  \,$
[repr\'esent\'e en trait fort et \`a l'aide des  zones gris\'ees]
pour les \'equations d'Euler-Saint Venant de la dynamique des gaz dans l'approximation isentropique.
Cas d'une entr\'ee supersonique.}
%% \bigskip
%%%%%%%%%%%%%%%%%%%%%%%%%%%%%%%%%%%%%%%%%%%%%%%%%%%%%%%%%%%%%%%%%%%%%%%%%%%%%%%%%%%

%%%%%%%%%%%%%%%%%%%%%%%%%%%%%%%%      figure 2-2    %%%%%%%%%%%%%%%%%%%%%%%%%%%%%%%% 
%%%%%%   \bigskip %% \centerline {  \epsfysize=9,0cm  \epsfbox  {../cnslim/Cnslim.fig2.epsf} }
%%%%%%   \centerline  {\includegraphics[width=.55\textwidth]   {fig-02.pdf}}
%%%%%%   \centerline { {\bf Figure 2}	\quad 	 Ensemble limite  $\,\, {\cal V} \bigl( W_0
%%%%%%   \bigr)  \,\,$  [repr\'esent\'e en trait fort }
%%%%%%   \centerline { et \`a l'aide des  zones gris\'ees] pour les \'equations d'Euler-Saint
%%%%%%   Venant }
%%%%%%   \centerline { de la dynamique des gaz dans l'approximation isentropique. }
%%%%%%   \centerline {  Cas d'une sortie subsonique ; le cas d'une entr\'ee subsonique est
%%%%%%   analogue.} 
%%%%%%%%%%%%%%%%%%%%%%%%%%%%%%%%%%%%%%%%%%%%%%%%%%%%%%%%%%%%%%%%%%%%%%%%%%%%%%%%%% 

%%%%%%%%%%%%%%%%%%%%%%%%%%%%%%%%%%%%%%%%%%%%%%%%%%%%%%%%%%%%%%%%%%%%%%%%%%%%%%%%%%% figure  2
%% \bigskip
\centerline  {\includegraphics[width=.55\textwidth]   {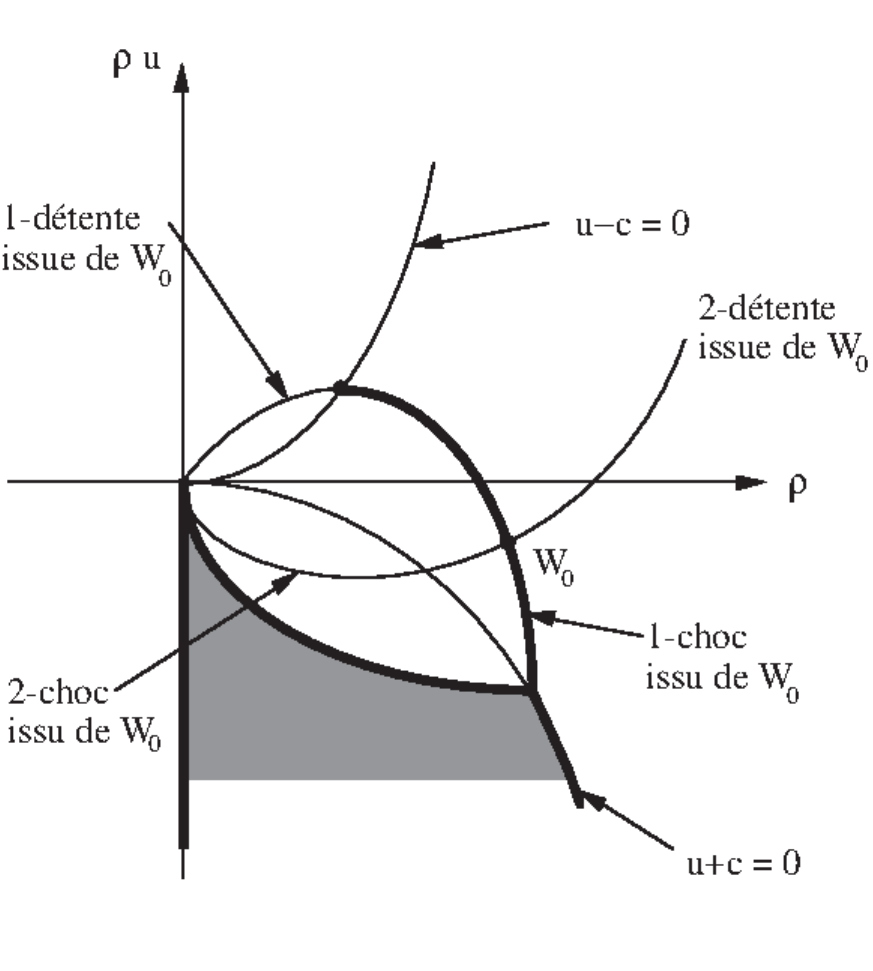}}
%%  \smallskip  
\vskip -.4 cm 
\noindent {{\bf Figure 2}. Ensemble limite  $\, {\cal V} \bigl( W_0 \bigr)  \,$
[repr\'esent\'e en trait fort et \`a l'aide des  zones gris\'ees]
pour les \'equations d'Euler-Saint Venant de la dynamique des gaz dans l'approximation isentropique.
Cas d'une sortie subsonique~; le cas d'une entr\'ee subsonique est
analogue.}  
\bigskip
%%%%%%%%%%%%%%%%%%%%%%%%%%%%%%%%%%%%%%%%%%%%%%%%%%%%%%%%%%%%%%%%%%%%%%%%%%%%%%%%%%%

\newpage

%%%%%%%%%%%%%%%%%%%%%%%%%%%%%%%%      figure 2-3    %%%%%%%%%%%%%%%%%%%%%%%%%%%%%%%% 
% \null\vskip -.2 cm
%%%  \centerline {  \epsfysize=9,0cm  \epsfbox  {../cnslim/Cnslim.fig3.epsf} } 
%%%%%%  \null\vskip -1.2  cm 
%%%%%%  \centerline { {\bf Figure 3}	\quad 	 Ensemble limite  $\,\, {\cal V} \bigl( W_0
%%%%%%  \bigr)  \,\,$  [repr\'esent\'e en trait fort }
%%%%%%  \centerline { et \`a l'aide des  zones gris\'ees] pour les \'equations d'Euler-Saint
%%%%%%  Venant }
%%%%%%  \centerline { de la dynamique des gaz dans l'approximation isentropique. }
%%%%%%  \centerline {  Cas d'une sortie supersonique de nombre de Mach inf\'erieur \`a $ \,
%%%%%%  {{2}\over{\gamma - 1}}.\, $ }
%%%%%%%%%%%%%%%%%%%%%%%%%%%%%%%%%%%%%%%%%%%%%%%%%%%%%%%%%%%%%%%%%%%%%%%%%%%%%%%%%% 

%%%%%%%%%%%%%%%%%%%%%%%%%%%%%%%%%%%%%%%%%%%%%%%%%%%%%%%%%%%%%%%%%%%%%%%%%%%%%%%%%%% figure  3
\bigskip
\centerline  {\includegraphics[width=.55\textwidth]   {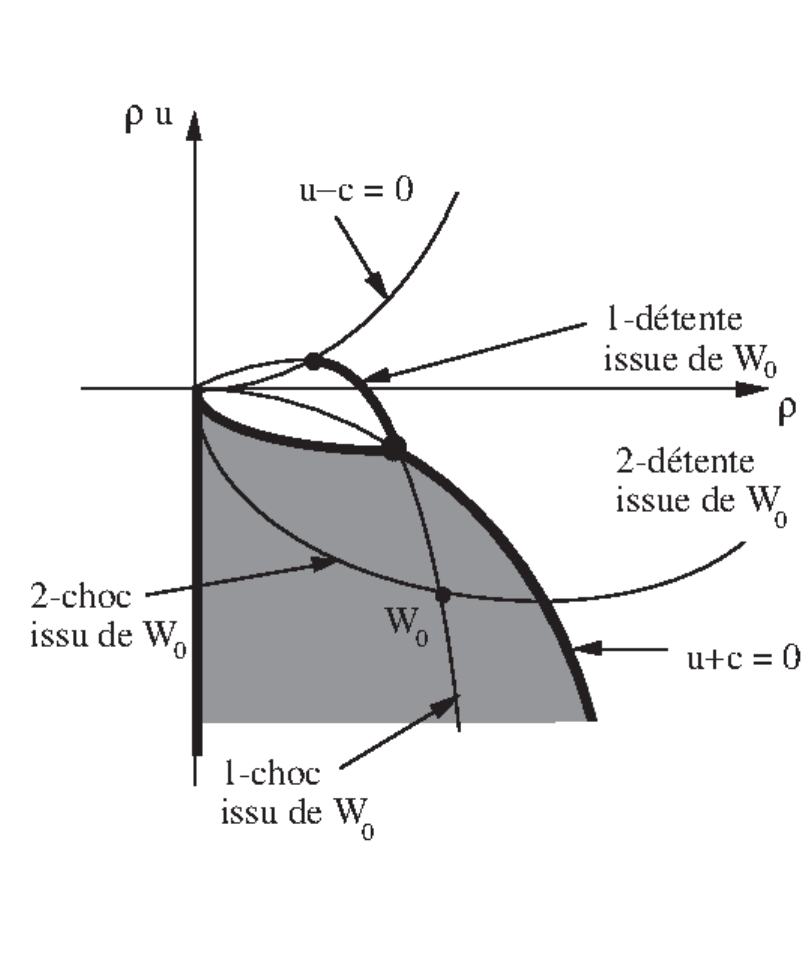}}
%% \smallskip  
\vskip -.8 cm 
\noindent {{\bf Figure 3}. Ensemble limite  $\, {\cal V} \bigl( W_0 \bigr)  \,$
[repr\'esent\'e en trait fort et \`a l'aide des  zones gris\'ees]
pour les \'equations d'Euler-Saint Venant de la dynamique des gaz dans l'approximation isentropique.
Cas d'une sortie supersonique de nombre de Mach inf\'erieur \`a $ \,{{2}\over{\gamma - 1}}$.}  
\bigskip
%%%%%%%%%%%%%%%%%%%%%%%%%%%%%%%%%%%%%%%%%%%%%%%%%%%%%%%%%%%%%%%%%%%%%%%%%%%%%%%%%%%

\bigskip \noindent   $\bullet \qquad  \,\,\, $  
Dans le cas d'un syst\`eme g\'en\'eral de lois de conservation, l'explicitation de
l'ensemble limite \'echoue car il est tr\`es difficile de manipuler pratiquement
toutes les entropies $\, \eta .\,$ Notons toutefois que pour le $p$-syst\`eme
(d\'efini par exemple dans Lax [La73]), Benabdallah-Serre [BS87] ont obtenu un
r\'egionnement non trivial de l'ensemble $\,  {\cal E} \bigl( w_0 \bigr) .\,$

\bigskip \noindent   $\bullet \qquad  \,\,\, $  
Nous avons par ailleurs remarqu\'e que dans les cas particuliers d'une \'equation
{\bf scalaire} (non n\'ecessairement convexe) et d'un {\bf syst\`eme hyperbolique
lin\'eaire}, l'ensemble limite $\,  {\cal E} \bigl( w_0 \bigr) \,$ est
caract\'eris\'e simplement \`a l'aide du {\bf pro\-bl\`eme de Riemann} ({\it c.f.} 
partie 1 et [DL87])~: 

\setbox11=\hbox{$\!\!$ valeurs en $\,\, {{x}\over{t}}=0^+ \,\,$ de la solution entropique}
\setbox12=\hbox {$\!\!$ du probl\`eme de Riemann $\,R(w_0 ,\, w) ,\, w \,$ variant}
\setbox40= \vbox {\halign{#&# \cr \box11 \cr \box12 \cr}}
\setbox41= \hbox{ $\vcenter {\box40} $}
\setbox44=\hbox{\noindent  (2.39) $\,\,\, \displaystyle  {\cal E} \bigl( w_0 \bigr) 
=  {\cal V} \bigl( w_0 \bigr)  \equiv  \left\{ \box41  \,\, \right\}  . $}  
\smallskip \noindent $ \box44 $

\smallskip \noindent
Une seconde formulation de la condition limite s'\'ecrit alors simplement 
[DL87]~: 

\smallskip \noindent  (2.40) $\qquad \displaystyle
w \bigl( t,\, 0^+ \bigr) \in {\cal V} \bigl( w_0(t) \bigr)  \qquad \qquad t > 0
\,.\,$ 

\smallskip \noindent 
Elle a l'avantage d'introduire un ensemble limite $\,\, {\cal V} \bigl( w_0 \bigr) 
\,\,$ {\bf calculable explicitement}. Remarquons que Benabdallah-Serre [BS87] ont
montr\'e que l'inclusion 

\smallskip \noindent  (2.41) $\qquad \displaystyle
{\cal V} \bigl( w_0 \bigr)  \, \, \subset    \, \, {\cal E} \bigl( w_0 \bigr) \,$ 

\smallskip \noindent
est toujours satisfaite, mais que l'\'egalit\'e (2.39) peut \^etre en d\'efaut. Par
ailleurs, Dubroca-Gallice [DG88] ont montr\'e que la m\'ethode de Glimm converge
lorsqu'on l'associe \`a la\br
premi\`ere condition limite (2.34) et que pour le
$p$-syst\`eme, la condition plus restrictive (2.40) conduit \'egalement \`a un
probl\`eme bien pos\'e.

\bigskip \noindent   $\bullet \qquad  \,\,\, $  
Avec Philippe Le Floch [DL87], nous avons calcul\'e et repr\'esent\'e graphi\-quement
l'ensemble $\,\, {\cal V} \bigl( w_0 \bigr)  \,\,$ dans le cas des \'equations
d'Euler-Saint Venant de la dynamique des gaz {\bf isentropiques} ainsi qu'illustr\'e
aux figures 1 \`a 3.   M\^eme lorsque l'\'etat $\, w_0 \,$ correspond \`a une entr\'ee
supersonique (figure 1), l'ensemble admissible $\,\, {\cal V} \bigl( w_0 \bigr)  \,\,$
n'est r\'eduit \`a $\, \{w_0\} \,$ que dans un {\bf voisinage} de $\, w_0 ,\,$ ce qui
correspond \`a l'approche lin\'earis\'ee. L'\'etat $\, w_0 \,$ peut \^etre reli\'e
\`a l'\'etat $\,\, w(t,\,0^+) \,\,$ par un 1-choc, ou bien $\,\, w(t,\,0^+) \,\,$ peut
correspondre \`a une sortie supersonique~(!).

\bigskip \noindent   $\bullet \qquad  \,\,\, $  
L'utilisation du probl\`eme de Riemann pour formuler la condition limite conduit \`a
un probl\`eme bien pos\'e lorsque les donn\'ees $\, w_0 \,$ et $\, v_0 \,$ sont des
\'etats constants~:  

\setbox11=\hbox{$\displaystyle \,\, \,\, {{\partial w}\over{\partial t}} \,+ \, 
{{\partial}\over{\partial x}} \, f(w) \,\, $  }
\setbox21=\hbox {$ = \,\, 0 \,$  }
\setbox31=\hbox {$ \qquad \qquad x > 0 \,,\qquad t > 0  \,$}
\setbox12=\hbox {$\displaystyle \,\, w(0,\,x) \,\,$  }
\setbox22=\hbox {$ = \,\, v_0  \, $}
\setbox32=\hbox {$ \qquad \qquad x > 0   \,$}
\setbox13=\hbox {$\displaystyle \,\,w(t,\,0^+) \,$  }
\setbox23=\hbox {$ \in \, {\cal V} \bigl( w_0 \bigr) \, $}
\setbox33=\hbox {$ \qquad \qquad t > 0  \,.  \,$}
\setbox40= \vbox {\halign{#&#&# \cr \box11 & \box21 & \box31 \cr \box12 & \box22  &
\box32 \cr \box13 & \box23  & \box33 \cr}}
\setbox41= \hbox{ $\vcenter {\box40} $}
\setbox44=\hbox{\noindent  (2.42) $\qquad \displaystyle \left\{ \box41 \right. $}  
\smallskip \noindent $ \box44 $

\smallskip  \smallskip \noindent
Nous pouvons par ailleurs envisager de fortes non lin\'earit\'es et pour les
\'equations d'Euler, seuls les probl\`emes d'apparition du vide limitent la
construction de la solution du probl\`eme de Riemann $\, R(w_g,\,w_d) .\,$ De plus,
le calcul des points fronti\`eres utilise alors les m\^emes outils que le calcul des
points int\'erieurs, comme nous l'\'etudions dans la seconde partie.

\bigskip \bigskip 
%%%%%%%%%%%%%%%%%%%%%%%%%%%%%%%%%%%%%%%%%%%%%%%%%%%%%%%%%%%%%%%%%%%%%%%%%%%%%%%  section 3
\noindent {\bf \large    3) \quad  Discr\'etisation des conditions aux limites} 
%%%%%%%%%%%%%%%%%%%%%%%%%%%%%%%%%%%%%%%%%%%%%%%%%%%%%%%%%%%%%%%%%%%%%%%%%%%%%%%%%%%%%%%%%%

\smallskip \noindent {\bf 3.1)  \quad  	 Diff\'erences finies ou volumes finis ? }

\noindent   $\bullet \qquad  \,\,\, $  
Nous notons $\, \Delta x \,$ (respectivement $\, \Delta t $) le pas d'espace
(respectivement de temps), suppos\'e uniforme, en vue de l'\'ecriture sous forme
discr\`ete de l'\'equation (2.1). Les sch\'emas conservatifs sont fond\'es sur une
\'ecriture des \'equations sous forme int\'egrale~: 

\smallskip \noindent  (3.1) $\qquad \displaystyle
{{1}\over{\Delta t}} \, \Bigl( \, W_j^{n\!+\!1} \,-\, W_j^{n} \, \Bigr) \,\,+\,\, 
{{1}\over{\Delta x}} \, \Bigl( \, f_{j\!+\!1/2} \,-\,  f_{j\!-\!1/2}  \, \Bigr)
\,\,=\,\, 0 \,$ 

\smallskip \noindent 
et pour les points int\'erieurs au domaine de calcul, $\,  W_j^{n} \,$ repr\'esente
une approximation de la valeur ponctuelle de $\, W \,$ au temps $\,\, t^n = n \, \Delta
t \,\,$ et au point $\, x_j \,$ aussi bien que la valeur moyenne du champ $\, W(t^n
,{\scriptstyle \bullet}) \,$ sur un intervalle de mesure $\, \,h = \Delta x \,\,$
autour du point $\, x_j .\,$ On peut raisonnablement dire que les deux m\'ethodes
co\"\i ncident  dans ce cas.

\bigskip \noindent   $\bullet \qquad  \,\,\, $  
Lorsqu'on se place sur un intervalle born\'e (figure 4), les deux m\'ethodes
conduisent \`a des probl\`emes sp\'ecifiques. Avec la m\'ethode des diff\'erences
finies, on construit en g\'en\'eral un ``sch\'ema fronti\`ere'' {\bf diff\'erent} du
sch\'ema (3.1), alors que la m\'ethode des volumes finis traite les cellules du bord
comme les autres et ne demande qu'une \'evaluation du {\bf flux fronti\`ere}. Nous
d\'etaillons ces deux approches dans les paragraphes qui suivent, en \'etudiant
successivement le cas des \'equations d'Euler {\bf lin\'earis\'ees}, le cas d'une
paroi et enfin celui d'une fronti\`ere fluide. La pr\'ecision des m\'ethodes
propos\'ees est variable selon les probl\`emes et les m\'ethodes utilis\'ees. La
m\'ethode des diff\'erences finies permet l'\'ecriture de sch\'emas au besoin tr\`es
pr\'ecis alors que celle des volumes finis d\'eg\'en\`ere au premier ordre dans la
plupart des exemples que nous proposons.

%%%%%%%%%%%%%%%%%%%%%%%%%%%%%%%%      figure 2-4    %%%%%%%%%%%%%%%%%%%%%%%%%%%%%%%%  
%%%%%%  \bigskip 
%%%  \centerline {  \epsfysize=3,0cm  \epsfbox  {../cnslim/Cnslim.fig4.epsf} } 
%%%%%% \smallskip  \smallskip
%%%%%% \centerline { {\bf Figure 4}	\quad 	 	Points de calcul pour la r\'esolution des
%%%%%% \'equations d'Euler }
%%%%%% \centerline { par la  m\'ethode des diff\'erences finies (en haut) :  $\,\,   x_j \, =
%%%%%% \,  j \, h \,\,$ }
%%%%%% \centerline {   et celle des volumes finis (en bas) :   $\,\,   x_j \, = \,  (j+1/2)\,
%%%%%% h . \,\,$ }
%%%%%%%%%%%%%%%%%%%%%%%%%%%%%%%%%%%%%%%%%%%%%%%%%%%%%%%%%%%%%%%%%%%%%%%%%%%%%%%%%%%%%  
 
%%%%%%%%%%%%%%%%%%%%%%%%%%%%%%%%%%%%%%%%%%%%%%%%%%%%%%%%%%%%%%%%%%%%%%%%%%%%%%%%%%% figure  4
\bigskip
\centerline  {\includegraphics[width=.40\textwidth]   {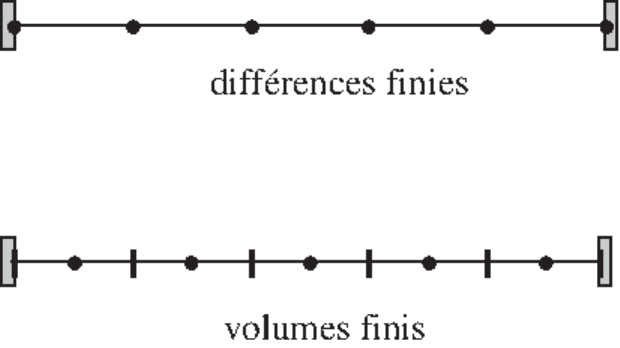}}
\smallskip  
\noindent {{\bf Figure 4}. Points de calcul pour la r\'esolution des \'equations d'Euler
par la  m\'ethode des diff\'erences finies (en haut)~:  $\,   x_j \, = \,  j \, h \,\,$ 
et celle des volumes finis (en bas)~:   $\,\,   x_j \, = \,  (j+1/2)\,  h $.}
\bigskip
%%%%%%%%%%%%%%%%%%%%%%%%%%%%%%%%%%%%%%%%%%%%%%%%%%%%%%%%%%%%%%%%%%%%%%%%%%%%%%%%%%%

\bigskip \noindent {\bf 3.2)  \quad  	 Equations d'Euler lin\'earis\'ees}

\noindent   $\bullet \qquad  \,\,\, $  
Nous nous limitons \`a la fronti\`ere de gauche (situ\'ee en $\, x = 0 $) d'un domaine
de calcul $\, \{ x > 0 \} .\,$ Une condition du type (2.26), {\it i.e.} 

\smallskip \noindent  (3.2) $\qquad \displaystyle
\varphi^{I} \,\,= \,\, \Sigma \, \varphi^{II} \,+\, g \qquad \qquad t \geq 0 \,,\,
\quad x = 0 \,$ 

\smallskip \noindent 
entre les variables caract\'eristiques entrantes $\,\,  \varphi^I \,\,$ et sortantes
$\,\, \varphi^{II} \,\,$ assure $\,p \,$ relations \`a la fronti\`ere (o\`u $p$ est
le nombre de valeurs propres positives de (2.24), {\it i.e.} le nombre de composantes de
$\,  \varphi^{I} $).

\bigskip \noindent   $\bullet \qquad  \,\,\, $  
La m\'ethode des diff\'erences finies demande donc $\, (n-p) \,$ relations
suppl\'e-mentaires qui constituent les ``conditions aux limites num\'eriques'' \`a
ajouter aux relations (3.2) pour calculer l'\'etat au bord. Diverses m\'ethodes
classiques (extrapolation en espace, en espace-emps, sch\'ema d\'ecentr\'e au bord,
etc...) sont d'utilisation courante (voir par exemple Yee, Beam et Warming [YBW82] ou
Cambier, Escande et Veuillot [CEV86]) et conduisent \`a des sch\'emas bien pos\'es au
sens de la stabilit\'e GKS (voir par exemple Goldberg-Tadmor [GT87] et les
r\'ef\'erences incluses). Nous d\'etaillons ici le proc\'ed\'e qui nous semble le
plus satis\-faisant, celui des {\bf relations de compatibilit\'e} (Viviand-Veuillot
[VV78], voir aussi Bramley-Sloan [BS77] et Kentzer [Ke71]). 

%%%%%%%%%%%%%%%%%%%%%%%%%%%%%%%%      figure 2-5    %%%%%%%%%%%%%%%%%%%%%%%%%%%%%%%%   
%%%  \bigskip \centerline {  \epsfysize=3,3cm  \epsfbox  {../cnslim/Cnslim.fig5.epsf} } 
%%%  \smallskip  \smallskip
%%%  \centerline { {\bf Figure 5}	\quad 	 Directions caract\'eristiques au bord du domaine
%%%  de calcul. }
%%%%%%%%%%%%%%%%%%%%%%%%%%%%%%%%%%%%%%%%%%%%%%%%%%%%%%%%%%%%%%%%%%%%%%%%%%%%%%%%%%   

%%%%%%%%%%%%%%%%%%%%%%%%%%%%%%%%%%%%%%%%%%%%%%%%%%%%%%%%%%%%%%%%%%%%%%%%%%%%%%%%%%% figure  5
\bigskip
\centerline  {\includegraphics[width=.60\textwidth]   {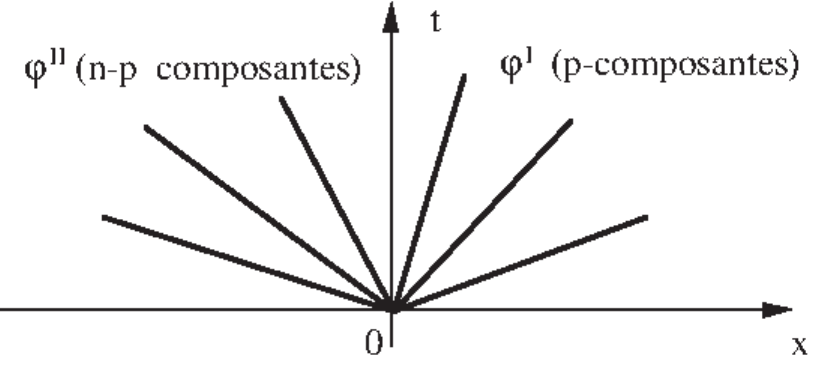}}
\smallskip  
\noindent {{\bf Figure 5}.  Directions caract\'eristiques au bord du domaine
de calcul.} 
\bigskip
%%%%%%%%%%%%%%%%%%%%%%%%%%%%%%%%%%%%%%%%%%%%%%%%%%%%%%%%%%%%%%%%%%%%%%%%%%%%%%%%%%%

\bigskip \noindent   $\bullet \qquad  \,\,\, $  
Les relations de compatibilit\'e reviennent, dans le cas lin\'eaire, \`a
discr\'etiser les \'equations du syst\`eme (2.24) qui correspondent aux {\bf ondes
sortantes} (composantes $\, \varphi^{II} )\,$~; on obtient ainsi $\, n-p \,$
\'equations d'advection~: 

\smallskip \noindent  (3.3) $\qquad \displaystyle
{{\partial \varphi^{II}}\over{\partial t}} \,+\, \Lambda^{II}({\overline W}) \, 
{{\partial \varphi^{II}}\over{\partial x}} \,\,= \,\, P^{II} \, C({\overline W},\,
\varphi) \,$ 

\smallskip \noindent 
o\`u $\,\,  P^{II} \,\,$ est le projecteur associ\'e aux variables $\, \,
\varphi^{II} \,\,$ dans la d\'ecomposition (2.25) (figure 5). Le sch\'ema final se
d\'eduit de (3.2)(3.3)~: une discr\'etisation {\bf d\'ecentr\'ee} de la d\'eriv\'ee
spatiale $\,\, {{\partial}\over{\partial x}} \,\,$  pr\'esente au sein de la relation
(3.3) permet de calculer $\,\,  \varphi^{II} \,\,$ en $\, x=0 \,$ au temps $\,
t^{n\!+\!1}  \,$ \`a partir de $\,\, \varphi \,\,$ au temps $ \, t^n \,$ et la
relation (3.2) ach\`eve la construction. L'\'evaluation des d\'eriv\'ees spatiales
\`a l'instant $\, t^{n\!+\!1}  \,$ conduit \`a divers sch\'emas implicites \`a la
fronti\`ere (voir par exemple Yee, Beam et Warming [YBW82]).

\bigskip \noindent   $\bullet \qquad  \,\,\, $  
La m\'ethode des volumes finis s'applique sans modification essentielle au bord de
$\, \Omega \,$~; il suffit de prendre $\, j=1 \,$ dans la relation (3.1) qui
d\'efinit le sch\'ema. Le flux aux interfaces est \'evalu\'e de fa\c{c}on
{\bf d\'ecentr\'ee} par une m\'ethode d'ordre un de Godunov [Go59] ou d'ordre deux de
Van Leer [VL79]. A la paroi, le flux $\,\,f_{j\!+\!1/2} \,\,$ est calcul\'e de
fa\c{c}on \`a prendre en compte la relation (3.2)~: 

\smallskip \noindent  (3.4) $\qquad \displaystyle
f_{1/2} \,\,= \,\, \Lambda^+({\overline W}) \, \bigl( \, \Sigma \, \varphi^{II}_1 \,+\,
g \, \bigr) \,\,+ \,\,  \Lambda^-({\overline W}) \,  \varphi^{II}_1 \,$ 

\smallskip \noindent
o\`u $\,\,  \varphi^{II}_1 \, \,$ d\'esigne la composante num\'ero $\, II \,$ de
l'\'etat $\, W_1 \,$ ou bien une valeur extrapol\'ee en $\, x=0^+ \,$ du champ \`a
partir des valeurs dans les premi\`eres mailles (pour les calculs d'ordre deux).
Remarquons que l'\'ecriture (3.4) de la condition limite (3.2) {\bf affaiblit} cette
derni\`ere puisque seul le flux num\'erique utilise la condition \`a la limite. Notons
\'egalement que l'\'etat paroi qui permet d'\'evaluer le flux gr\^ace \`a la relation
(3.4) r\'esulte d'une extrapolation en espace des variables sortantes $\,\,
\varphi^{II} \,\,$ et d'un calcul des variables entrantes $\,\, \varphi^I \,\,$ \`a
partir de la relation  (3.2), avant l'incr\'ementationen temps (3.1), alors que
l'approche par la m\'ethode des diff\'erences finies jointe aux relations de
compatibilit\'e revient \`a coupler les deux \'etapes pour les composantes sortantes,
avant de recalculer $\, \, \varphi^{II} \,\,$ gr\^ace \`a la condition limite (3.2).
L'avantage de l'approche ``volumes finis'' pour les calculs fronti\`eres tient au fait 
qu'aucun sch\'ema num\'erique suppl\'ementaire n'est n\'ecessaire (pour les calculs
au premier ordre en espace au moins~!) pour prendre en compte la fronti\`ere, comme
l'avaient remarqu\'e Godunov   {\it et al.}    [GZIKP79] et Rizzi [Ri81] par exemple. Signalons
aussi que les \'etudes  th\'eoriques de stabilit\'e \`a la fronti\`ere sont beaucoup
moins d\'evelopp\'ees avec cette approche qu'avec l'approche ``diff\'erences finies''
classique. 

\bigskip \noindent {\bf 3.3)  \quad  	 Fronti\`ere fluide}

\noindent   $\bullet \qquad  \,\,\, $  
Comme on l'a vu dans la premi\`ere partie, on distingue classiquement quatre cas
selon le type d'\'ecoulement \`a la fronti\`ere, c'est \`a dire le nombre $p$ de
composantes de $\, \varphi^I .\, $ L'approche aux diff\'erences finies remplace la
relation (3.2) par $\,p \,$ relations {\bf non lin\'eaires} adapt\'ees au probl\`eme, 
{\it i.e.} \'etat impos\'e,  enthalpie totale et entropie impos\'ees,  pression impos\'ee
ou   pas de relation pour les quatre cas classiques~: 

\smallskip \noindent  (3.5) $\qquad \displaystyle
B(W) \,\,= \,\, 0 \qquad \qquad \,(p \,\,$ relations non lin\'eaires). 

\smallskip \noindent
Les $\,(n-p) \,$ relations suppl\'ementaires sont le plus souvent issues des {\bf
relations de compatibilit\'e} (Viviand-Veuillot [VV78])~: le syst\`eme (2.1) admet
des vecteurs propres $\, r_{j}(W) \,$ (2.10) et des formes lin\'eaires propres
$\,\, \ell_{j}(W) \, \,$ (ou ``vecteurs propres \`a gauche''~; $\,
\ell_{j}(W)^{\displaystyle \rm t} \,$ est en fait vecteur propre de $\, A(W) 
^{\displaystyle \rm t} ) , \,$  d\'efinies par les relations 

\smallskip \noindent  (3.6) $\qquad \displaystyle
\ell_{j}(W) \,  \smb  \, A(W) \,\,= \,\, \lambda_{j}(W) \,\, \ell_{j}(W) \,. \,$ 

\smallskip \noindent 
Le syst\`eme (2.1) est r\'e\'ecrit sous la forme non conservative \'equivalente 

\smallskip \noindent  (3.7) $\qquad \displaystyle
\ell_{j}(W) \, \smb \,{{\partial W}\over{\partial t}} \,+\, 
\lambda_{j}(W) \, \ell_{j}(W) \, \smb \, {{\partial W}\over
{\partial x}} \,\,= \,\, 0 \,\,, \qquad j \,=\, 1,\, \cdots \,, \, n \, $ 

\smallskip \noindent 
qui constitue l'ensemble des relations de compatibilit\'e sous forme non lin\'eaire.
On ne conserve pour d\'efinir le sch\'ema limite que les indices $\,j \,$ qui
correspondent \`a des directions caract\'eristiques sortantes $\, (j=1,\, \cdots ,\,
n-p) \,$ et l'on a en d\'efinitive~: 

\setbox11=\hbox {$\displaystyle \,\, B_{j}(W) \,\,= \,\, 0 \,   $}
\setbox21=\hbox {$  \qquad j \,= \, 1 \,,\, \cdots \,,\, p \,  $}
\setbox12=\hbox {$\displaystyle \,\, \ell_{j}(W) \, \smb \,
{{\partial W}\over{\partial t}} \,+\,  \lambda_{j}(W) \, \ell_{j}(W) \,{\scriptstyle
\bullet}\, {{\partial W}\over {\partial x}} \,\,= \,\, 0\,  $}
\setbox22=\hbox {$  \qquad  j \,= \, 1 \,,\, \cdots \,,\, n-p \,. \,  $}
\setbox40= \vbox {\halign{#&# \cr \box11 & \box21\cr \box12 & \box22 \cr  \cr}}
\setbox41= \hbox{ $\vcenter {\box40} $}
\setbox44=\hbox{\noindent  (3.8) $\quad \displaystyle \left\{ \box41 \right. $}  
\smallskip \noindent $ \box44 $

\smallskip \noindent 
Lorsqu'on utilise la valeur de grille $\,\, W_0^n \,$ (sur la fronti\`ere du domaine
de calcul) dans l'expression de la forme lin\'eaire $\, \ell_{j} \,$ [on approche $\, 
\ell_{j}(W) \,$ par $\,\, \ell_{j}(W_0^n) $], les relations (3.7) sont identiques
aux relations (3.3) lin\'earis\'ees et fournissent les sch\'emas explicites les plus
utilis\'es (Viviand-Veuillot [VV78], Cambier-Escande-Veuillot [CEV86] et les
r\'ef\'erences incluses). Lorsqu'on lin\'earise le terme en $\,
{{\partial}\over{\partial t}} \,$ de (3.8) autour de la valeur de grille $\,  W_0^n
,\,$ (Chakravarthy [Ch83]), les relations obtenues d\'efinissent un sch\'ema
{\bf implicite} pour l'incr\'ement $\,\, \bigl(  W_0^{n\!+\!1} -  W_0^n \bigr) \,\,$
du point fronti\`ere.

%%%%%%%%%%%%%%%%%%%%%%%%%%%%%%%%      figure 2-6    %%%%%%%%%%%%%%%%%%%%%%%%%%%%%%%%    
%%%  \bigskip \centerline {  \epsfysize=3,5cm  \epsfbox  {../cnslim/Cnslim.fig6.epsf} } 
%%% \smallskip  \smallskip
%%%%%%  % fin de la version linux
%%%%%%  \centerline { {\bf Figure 6}	\quad 		Sortie subsonique faiblement non lin\'eaire }
%%%%%%  \centerline { (d'apr\`es Osher-Chakravarthy [OC83]). }
%%%%%%%%%%%%%%%%%%%%%%%%%%%%%%%%%%%%%%%%%%%%%%%%%%%%%%%%%%%%%%%%%%%%%%%%%%%%%%%%%%    

%%%%%%%%%%%%%%%%%%%%%%%%%%%%%%%%%%%%%%%%%%%%%%%%%%%%%%%%%%%%%%%%%%%%%%%%%%%%%%%%%%% figure  6
\bigskip
\centerline  {\includegraphics[width=.64\textwidth]   {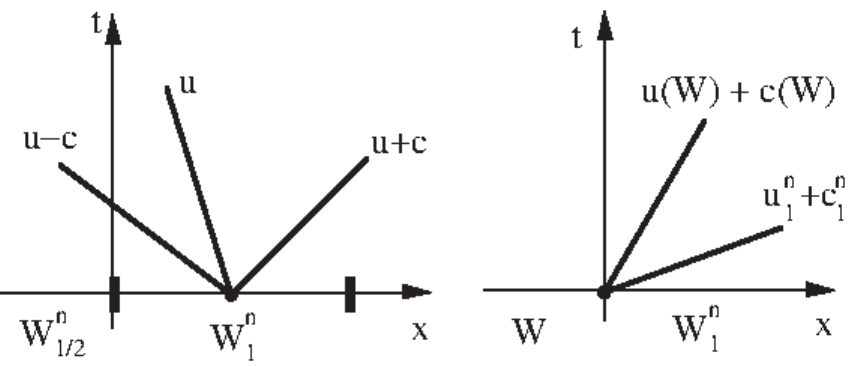}}
\smallskip  
\noindent {{\bf Figure 6}. Sortie subsonique faiblement non lin\'eaire
(d'apr\`es Osher-Chakravarthy [OC83]).}
%% \bigskip
%%%%%%%%%%%%%%%%%%%%%%%%%%%%%%%%%%%%%%%%%%%%%%%%%%%%%%%%%%%%%%%%%%%%%%%%%%%%%%%%%%%

\bigskip \noindent   $\bullet \qquad  \,\,\, $  
Pour les m\'ethodes de volumes finis, Osher-Chakravarthy [OC83] ont propos\'e un
proc\'ed\'e que nous g\'en\'eraliserons dans la troisi\`eme partie. L'id\'ee consis\-te
\`a calculer l'\'etat fronti\`ere $\, \, W_{1/2} \,\,$ par une \'etude des ondes
non lin\'eaires d'un probl\`eme de Riemann et non plus par les relations de
compatibilit\'e (3.7). L'\'etat fronti\`ere $\,  \, W_{1/2} \,\,$ v\'erifie
toujours la condition limite (3.5) et est reli\'e \`a l'\'etat $\,\, W_1^n \,\,$ dans
la premi\`ere maille par une famille de $p$ ondes simples ({\it c.f.} (2.17)(2.18)) et la
figure 6)~; il appartient donc \`a une vari\'et\'e de codimension $p$ (qui passe par
l'\'etat $\,  \, W_{1/2} $)  et satisfait aux $p$ \'equations (3.5), ce qui pose 
{\it a priori}    correctement le probl\`eme. Nous d\'etaillons le cas d'une sortie subsonique
o\`u la pression de sortie est suppos\'ee fix\'ee (voir \'egalement
Osher-Chakravarthy [OC83]) et la vari\'et\'e est alors de codimension $\, p = 1 .\,$
Nous utilisons le solveur du probl\`eme de Riemann propos\'e par Osher [Os81], qui
n'utilise que des ondes de d\'etentes, \'eventuellement multivalu\'ees (voir Van Leer
[VL84], Osher [Os84] et divers d\'etails dans  [Du87]). 

\bigskip \noindent   $\bullet \qquad  \,\,\, $  
On d\'etermine d'abord l'\'etat $\, W \,$ de pression impos\'ee $\,\, {\overline p}
\,\,$ tel que $\, W \,$ est li\'e \`a l'\'etat $\,\, W_1^n \,\,$ par une 3-onde de
d\'etente ({\it c.f.} 1.14))~: 

\setbox11=\hbox{$\displaystyle \,\,p(W) \,\, $  }
\setbox21=\hbox{$\displaystyle \,\, = \,\, {\overline p} \,$  }
\setbox12=\hbox {$\displaystyle \,\, u(W) \,-\, {{2}\over{\gamma \!-\!1}} \, c(W)
\,\, $}
\setbox22=\hbox{$\displaystyle \,\,= \,\, u_1^n \,-\,  {{2}\over{\gamma \!-\!1}} \,
c_1^n \, $}
\setbox13=\hbox {$\displaystyle \,\, {{p(W)} \over{\rho(W)^{\gamma}}} \,\, $}
\setbox23=\hbox{$\displaystyle \,\,= \,\, {{p_1^n} \over{{\rho_1^n}^{\gamma}}} \,.\,$}
\setbox40= \vbox {\halign{#&# \cr \box11 & \box21  \cr \box12   & \box22 \cr \box13 
& \box23  \cr}}
\setbox41= \hbox{ $\vcenter {\box40} $}
\setbox44=\hbox{\noindent  (3.9) $\qquad \displaystyle \left\{ \box41 \right. $}  
\smallskip \noindent $ \box44 $

\smallskip \noindent 
Si l'on suppose $\, u_1^n <  0 <   u_1^n + c_1^n $,  l'\'etat $ \, W_1^n \,$
dans la premi\`ere cellule satisfait aux in\'egalit\'es de sortie subsonique
(figure 6). Il s'agit donc d'un cas de {\bf faible non-lin\'earit\'e} et c'est
l'hypoth\`ese faite le plus souvent dans les applications (voir par exemple
Hemker-Spekreijse [HS86]). Dans ce cas, l'\'etat $\,\, W_{1/2} \,\,$ de paroi est
exactement l'\'etat interm\'ediaire calcul\'e en (3.9) et l'on a~: 

\smallskip \noindent  (3.10) $\qquad \displaystyle
f_{1/2} \,\, = \,\, f(W) \,$ 

\smallskip \noindent 
o\`u le flux $\, f({\scriptstyle \bullet})  \,$ est d\'etermin\'e \`a la relation
(2.3).

%%%%%%%%%%%%%%%%%%%%%%%%%%%%%%%%      figure 2-7    %%%%%%%%%%%%%%%%%%%%%%%%%%%%%%%%  
%%%   \bigskip \centerline {  \epsfysize=4,0cm  \epsfbox  {../cnslim/Cnslim.fig7.epsf} } 
%%%%%%%  \smallskip  \smallskip
% fin de la version linux
%%%%%%%  \centerline { {\bf Figure 7}	\quad 		Sortie subsonique fortement non lin\'eaire } 
%%%%%%%  %%%%%%%  \centerline {  calcul\'ee gr\^ace au sch\'ema d'Osher [Os81]. }
%%%%%%%%%%%%%%%%%%%%%%%%%%%%%%%%%%%%%%%%%%%%%%%%%%%%%%%%%%%%%%%%%%%%%%%%%%%%%%%%%%

%%%%%%%%%%%%%%%%%%%%%%%%%%%%%%%%%%%%%%%%%%%%%%%%%%%%%%%%%%%%%%%%%%%%%%%%%%%%%%%%%%% figure  7
\bigskip
\centerline  {\includegraphics[width=.40\textwidth]   {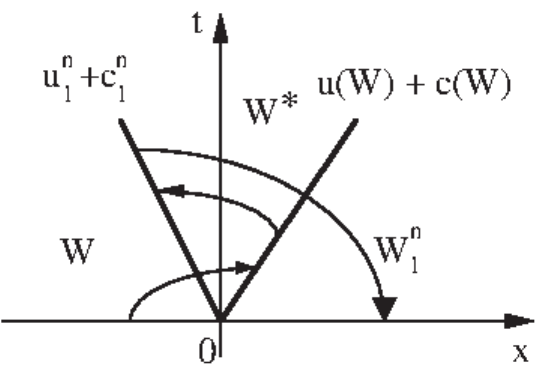}}
\smallskip  
\noindent {{\bf Figure 7}. Sortie subsonique fortement non lin\'eaire
calcul\'ee gr\^ace au sch\'ema d'Osher [Os81]. }
%% \bigskip
%%%%%%%%%%%%%%%%%%%%%%%%%%%%%%%%%%%%%%%%%%%%%%%%%%%%%%%%%%%%%%%%%%%%%%%%%%%%%%%%%%%

\bigskip \noindent   $\bullet \qquad  \,\,\, $  
Si, au contraire, on suppose  $\,\, u_1^n + c_1^n \, \leq \, 0 ,\,$ l'\'etat $\,\, 
W_1^n \,\,$ n'est pas lui-m\^eme un \'etat correspondant \`a une sortie subsonique. Il
y a donc {\bf forte non lin\'earit\'e} \`a la fronti\`ere et dans ce cas, les travaux
de Osher-Chakravarthy [OC83] donnent d\'ej\`a quelques \'el\'ements de r\'eponse.
L'\'etat $\, W \,$ est reli\'e \`a $\,\,  W_1^n \,\,$ par une {\bf d\'etente
\'eventuellement multivalu\'ee} et il faut abandonner la notion d'\'etat
fronti\`ere si l'on utilise le flux d'Osher. Il suffit de d\'eterminer le flux
\`a l'entr\'ee du domaine \`a l'aide de $\, W, \,W_1^n \,\,$ et de l'\'etat sonique
$\, W^* \,$ le long de la 3-d\'etente, calcul\'ee gr\^ace aux relations 

\setbox11=\hbox{$\displaystyle \,\,u(W^*) \,+\, c(W^*) \,\,  $  }
\setbox21=\hbox{$\displaystyle \,\, = \,\, 0 \,$  }
\setbox12=\hbox {$\displaystyle \,\, u(W^*) \,-\, {{2}\over{\gamma \!-\!1}} \, c(W^*)
\,\, $}
\setbox22=\hbox{$\displaystyle \,\,= \,\, u_1^n \,-\,  {{2}\over{\gamma \!-\!1}} \,
c_1^n \, $}
\setbox13=\hbox {$\displaystyle \,\, {{p(W^*)} \over{\rho(W^*)^{\gamma}}} \,\, $}
\setbox23=\hbox{$\displaystyle \,\,= \,\, {{p_1^n} \over{{\rho_1^n}^{\gamma}}} \,.\,$}
\setbox40= \vbox {\halign{#&# \cr \box11 & \box21  \cr \box12   & \box22 \cr \box13 
& \box23  \cr}}
\setbox41= \hbox{ $\vcenter {\box40} $}
\setbox44=\hbox{\noindent  (3.11) $\qquad \displaystyle \left\{ \box41 \right. $}  
\smallskip \noindent $ \box44 $

\smallskip \noindent
Avec Osher-Chakravarthy [OC83], nous supposons que l'\'etat $\, W \,$ v\'erifie la condition\br
$  u(W) + c(W) \geq 0 \,$~; la 3-onde de d\'etente ``contient'' la
fronti\`ere fluide $\, {{x}\over{t}} = 0 \,$ ({\it c.f.}  figu\-re 7) et le flux num\'erique
de fronti\`ere s'exprime selon 

\smallskip \noindent  (3.12) $\qquad \displaystyle
f_{1/2} \,\,= \,\, f(W) \,- \, f(W^*) \,+\, f(W_1^n ) \,. \,$ 

\smallskip \noindent 
D'autres cas de figure doivent \^etre envisag\'es pour prendre en compte les
diff\'e\-rentes configurations que peut prendre la 3-d\'etente (\'eventuellement
multivalu\'ee). Nous reviendrons dans la troisi\`eme partie sur une pr\'esentation
g\'en\'erale de ces probl\`emes. 

\bigskip \noindent   $\bullet \qquad  \,\,\, $  
Le cas d'une fronti\`ere fluide peut donc \^etre trait\'e dans de nombreux cas non
lin\'eaires comme une extension du cas lin\'eaire, o\`u les conditions aux limites
non-lin\'eaires (3.5) sont coupl\'ees aux relations de compatibilit\'e (3.7)
associ\'ees aux  caract\'eristiques sortantes. Cette approche est traditionnelle pour
les sch\'emas de volumes finis, le calcul de l'\'etat du bord propos\'e par
Osher-Chakravarthy [OC83], par r\'esolution d'un {\bf probl\`eme de Riemann partiel},
ne d\'efinit le sch\'ema \`a la fronti\`ere que par l'introduction du flux num\'erique
associ\'e. 

\bigskip \noindent   $\bullet \qquad  \,\,\, $  
Signalons enfin le probl\`eme des conditions aux limites {\bf absorbantes}, o\`u il
faut expri\-mer que les ondes quittent le domaine de calcul sans r\'eflexion. L'\'etude
monodimensionnelle de Hedstrom [He79] revient \`a \'ecrire que l'\'etat \`a la
fronti\`ere est une ``combinaison'' d'ondes de d\'etentes {\bf sortantes}, ce qui
implique (Hedstrom [He79])~: 

\smallskip \noindent  (3.13) $\qquad \displaystyle
\ell_{j} \, {\scriptstyle \bullet} \, {{\partial W}\over{\partial t}} \,\, = \,\, 0
\,.\,$ 

\smallskip \noindent 
Les relations (3.13) sont \`a substituer \`a (3.5) dans une formulation non
lin\'eaire de type (3.8) pour les sch\'emas aux diff\'erences finies. Pour une
extension bidimensionnelle, nous renvoyons \`a Thompson [Th87]. 

\bigskip \noindent {\bf 3.4)  \quad  	 Parois solides}

 \noindent   $\bullet \qquad  \,\,\, $  
Ce cas de condition limite est physiquement tr\`es diff\'erent du pr\'ec\'edent
puisqu'un obstacle est pr\'esent dans l'\'ecoulement et la condition limite physique

\smallskip \noindent  (3.14) $\qquad \displaystyle
u \,=\, 0 \qquad \qquad  (u \, $~: vitesse de l'\'etat paroi $ \, W $)

\smallskip \noindent
exprime la non-p\'en\'etrabilit\'e du fluide \`a la paroi.

%%%%%%%%%%%%%%%%%%%%%%%%%%%%%%%%      figure 2-8  %%%%%%%%%%%%%%%%%%%%%%%%%%%%%%%%  
%%%  \bigskip \centerline {  \epsfysize=1,5cm  \epsfbox  {../cnslim/Cnslim.fig8.epsf} } 
%%%%%  \smallskip  \smallskip
% fin de la version linux
%%%%%  \centerline { {\bf Figure 8}	\quad 	Point fictif dans la paroi (sch\'emas aux
%%%%%  diff\'erences finies). }
%%%%%%%%%%%%%%%%%%%%%%%%%%%%%%%%%%%%%%%%%%%%%%%%%%%%%%%%%%%%%%%%%%%%%%%%%%%%%%%%%%%  

%%%%%%%%%%%%%%%%%%%%%%%%%%%%%%%%%%%%%%%%%%%%%%%%%%%%%%%%%%%%%%%%%%%%%%%%%%%%%%%%%%% figure  8
\bigskip
\centerline  {\includegraphics[width=.45\textwidth]   {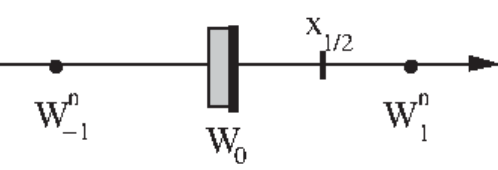}}
\smallskip  
\noindent {{\bf Figure 8}. Point fictif dans la paroi (sch\'emas aux
diff\'erences finies).}
\bigskip
%%%%%%%%%%%%%%%%%%%%%%%%%%%%%%%%%%%%%%%%%%%%%%%%%%%%%%%%%%%%%%%%%%%%%%%%%%%%%%%%%%%

\bigskip \noindent   $\bullet \qquad  \,\,\, $  
Le traitement de (3.13) \`a l'aide des sch\'emas aux diff\'erences a \'et\'e
propos\'e tr\`es t\^ot. L'approche traditionnelle propose d'introduire un point
fictif $\, W_{-1} \,$ dans la paroi (figure 8) et d'\'ecrire une condition de
sym\'etrie 

\setbox11=\hbox{$\displaystyle \,\,u_{-1} \,\,  $  }
\setbox21=\hbox{$\displaystyle \,\, = \,\, -u_1 \, $  }
\setbox12=\hbox {$\displaystyle \,\,p_{-1} \,\,  $  }
\setbox22=\hbox{$\displaystyle \,\,= \,\,p_1 \, $  }
\setbox13=\hbox {$\displaystyle \,\, \rho_{-1} \,\,  $  }
\setbox23=\hbox{$\displaystyle \,\,= \,\,\rho_{1} \,\,  $  }
\setbox40= \vbox {\halign{#&# \cr \box11 & \box21  \cr \box12   & \box22 \cr \box13 
& \box23  \cr}}
\setbox41= \hbox{ $\vcenter {\box40} $}
\setbox44=\hbox{\noindent  (3.15) $\qquad \displaystyle \left\{ \box41 \right. $}  
\smallskip \noindent $ \box44 $

\smallskip \noindent
de fa\c{c}on \`a incr\'ementer l'\'etat paroi $\, W_0 \,$ \`a l'aide du sch\'ema \`a
trois points utilis\'e pour les points int\'erieurs. Mais Moretti [Mo68] a montr\'e
que les conditions (3.15) introduisent des conditions suppl\'ementaires \`a la paroi
qui sur-sp\'ecifient le probl\`eme continu. Aussi Roache [Ro72] recommende-t-il
d'utiliser les volu\-mes finis (``second mesh system'').  L'approche courante actuellement
depuis Viviand-Veuillot [VV78] utilise la relation de compatibilit\'e correspondant
\`a la valeur propre $\, \lambda_{1}  \,$(sortant  du domaine de calcul) pour
calculer la pression paroi (essentiellement pour les applications a\'erodynamiques~!).
Pour les d\'etails, nous renvoyons \`a Cambier-Escande-Veuillot [CEV86]. 

\bigskip \noindent   $\bullet \qquad  \,\,\, $  
L'utilisation des volumes finis est plus facile \`a formuler~; compte tenu de la
condition de non-p\'en\'etrabilit\'e (3.14), le flux fronti\`ere $\,\, f_{1/2}
\,\,$ prend la forme alg\'ebrique suivante ({\it c.f.}  (2.3))~: 

\smallskip \noindent  (3.16) $\qquad \displaystyle
f_{1/2} \,\,= \,\, \bigl( \, 0 \,,\, p_{1/2} \,,\, 0 \, \bigr)^{\displaystyle \rm
t} \,$ 

\smallskip \noindent
et la pression paroi $\, p_{1/2} \,$ d\'efinit compl\`etement le sch\'ema. Dans ce
cas, la m\'ethode de l'``\'etat miroir'' (condition de sym\'etrie) d\'efinit encore un
\'etat fictif $\,\, W_0 \,\,$ de sorte que~: 

\setbox11=\hbox{$\displaystyle \,\,u_{0} \,\,  $  }
\setbox21=\hbox{$\displaystyle \,\, = \,\, -u_1 \, $  }
\setbox12=\hbox {$\displaystyle \,\,p_{0} \,\,  $  }
\setbox22=\hbox{$\displaystyle \,\,= \,\,p_1 \, $  }
\setbox13=\hbox {$\displaystyle \,\, \rho_{0} \,\,  $  }
\setbox23=\hbox{$\displaystyle \,\,= \,\,\rho_{1} \,. \,  $  }
\setbox40= \vbox {\halign{#&# \cr \box11 & \box21  \cr \box12   & \box22 \cr \box13 
& \box23  \cr}}
\setbox41= \hbox{ $\vcenter {\box40} $}
\setbox44=\hbox{\noindent  (3.17) $\qquad \displaystyle \left\{ \box41 \right. $}  
\smallskip \noindent $ \box44 $

\smallskip \noindent
Cet \'etat ne sert qu'\`a l'\'evaluation du flux paroi (3.16), par l'interm\'ediaire
d'une r\'esolution exacte du probl\`eme de Riemann $\,\, R(W_0,\, W_1) .\,$ En effet,
l'extrapolation de la pression et de la densit\'e propos\'es dans (3.17) ne sont pas
introduits explicitement dans le sch\'ema num\'erique, ce qui n'est pas le cas avec
l'approche aux diff\'erences finies (relations (3.15)). 

\bigskip \noindent   $\bullet \qquad  \,\,\, $  
Notons \'egalement que dans le cas de faibles non lin\'earit\'es \`a la paroi (le
plus courant~; la vitesse normale $\, u_1 \,$ dans la cellule jouxtant le bord est
``petite''), le calcul de la pression $\,\,p_{1/2} \,\,$ peut \^etre effectu\'e par
extrapolation \`a l'aide d'un sch\'ema aux diff\'erences pr\'ecis au second ordre
(voir par exemple Lerat [Le81]). Signalons enfin l'approche utilis\'ee \`a l'Inria
(Stoufflet [St84] par exemple). Celle-ci s'apparente \`a la fois aux diff\'erences
finies (ou aux \'el\'ements finis~!) puisque le n\oe ud \`a calculer est situ\'e sur
la paroi et aux volumes finis puisque le sch\'ema utilis\'e, du premier ordre en
espace, revient \`a \'ecrire un bilan dans la demi-maille $\,\, [x_0,\, x_{1/2} ] \,$
(figure 8). Le flux paroi est \'evalu\'e gr\^ace \`a la relation (3.16) et la
pression correspondante est simplement la pression de l'\'etat $\, W_0 .\,$

\bigskip \bigskip 
%%%%%%%%%%%%%%%%%%%%%%%%%%%%%%%%%%%%%%%%%%%%%%%%%%%%%%%%%%%%%%%%%%%%%%%%%%%%%%%  section 4
\noindent {\bf \large    4) \quad Probl\`eme de Riemann partiel \`a la fronti\`ere} 
%%%%%%%%%%%%%%%%%%%%%%%%%%%%%%%%%%%%%%%%%%%%%%%%%%%%%%%%%%%%%%%%%%%%%%%%%%%%%%%%%%%%%%%%%%

\noindent   $\bullet \qquad  \,\,\, $  
Nous pr\'esentons dans cette derni\`ere partie une technique g\'en\'erale qui permet
la prise en compte des fortes non lin\'earit\'es aux fronti\`eres du domaine de
calcul, lorsqu'on utilise la m\'ethode des volumes finis et les sch\'emas ``de type
Godunov'' (au sens de Harten-Lax et Van Leer [HLV83]). Dans le premier paragraphe,
nous exposons l'ensemble du sch\'ema dans le cas d'une pr\'ecision du {\bf premier
ordre}, puis nous montrons (au second paragraphe) comment le cas particulier de
faibles effets non lin\'eaires conduit aux r\'esultats classiques expos\'es plus
haut. Nous terminons par un cas test monodimensionnel o\`u de fortes non
lin\'earit\'es sont pr\'esentes. 

\bigskip \noindent {\bf 4.1)  \quad  	 Volumes finis prenant en compte la fronti\`ere}

\noindent   $\bullet \qquad  \,\,\, $  
Nous rappelons que nous cherchons une approximation $\,\,W_j^n \,\,$ de la valeur
moyenne des variables conservatives (2.2) dans la maille $\,\, K_j \,=\, ] \,
(j\!-\!1/2) \, \Delta x \,,\, (j\!+\!1/2) \, \Delta x \, [ \,\,$ au temps $\, t^n = n
\, \Delta t .\,$ Le sch\'ema de Godunov [Go59] ou les m\'ethodes de type Godunov
([Harten-Lax-Van Leer [HLV83]) consistent \`a int\'egrer la loi de conservation (2.1)
dans le domaine d'espace-temps $\,\,] n\, \Delta t \,,\, (n\!+\!1) \, \Delta t [ \,
\times \, K_j \,\,$ et on l'\'ecrit usuellement~: 

\smallskip \noindent  (4.1) $\qquad \displaystyle
{{1}\over{\Delta t}} \Bigl( \, W_j^{n\!+\!1} \,- \, W_j^n \, \Bigr) \,+\,
{{1}\over{\Delta x}} \Bigl( \,f_{j\!+\!1/2} \,-\, f_{j\!-\!1/2} \, \Bigr) \,\,\,=
\,\,\, 0 \qquad j \,=\, 1 \,,\, \cdots \,,\, N \,. \,$ 

\smallskip \noindent
Pour les mailles internes $\,\, \bigl(j\!+\!{1\over2} \,=\, {3\over2},\, {5\over2},\, 
\cdots ,\,  N-{1\over2} \, \bigr)  , \,\,$ le flux num\'erique correspondant est
calcul\'e par r\'esolution exacte (Godunov) ou approch\'ee (Roe [Ro81], Osher [Os81],
Collela-Glaz [CG85]) du probl\`eme de Riemann $\,\, R(W_j^n ,\, W_{j\!+\!1}^n) \,;$
nous notons $\, \Phi \,$ le flux num\'erique correspondant~: 

\smallskip \noindent  (4.2) $\qquad \displaystyle
f_{j\!+\!1/2} \,\,=\,\, \Phi \bigl( W_j^n ,\, W_{j\!+\!1}^n \bigr) \,, \qquad \qquad
j \,=\, 1 ,\, 2 ,\, \cdots ,\, N-1 \,. \,$ 

\smallskip \noindent 
Nous supposons que le flux num\'erique $\, \Phi \,$ est calcul\'e gr\^ace aux 
{\bf ondes} d'un probl\`eme de Riemann ({\it c.f.}  (2.17)(2.18)) et non par une
d\'ecomposition de flux, ce qui en pratique nous limite aux solveurs de Godunov
(exact), Osher (ondes de d\'etente seulement) ou Collela-Glaz (ondes de choc
seulement)~; le solveur de Roe (discontinuit\'es de contact seulement) ne d\'efinit
pas \`a notre connaissance de $j$-ondes $\,\,  U_{j} \, \,$ dans l'espace des \'etats
et nous ne l'utilisons donc pas dans la suite. 

\bigskip \noindent   $\bullet \qquad  \,\,\, $  
A la fronti\`ere du domaine de calcul, les conditions aux limites conduisent \`a se
donner certains param\`etres caract\'eristiques (voir les relations (3.5)). Nous
appelons {\bf vari\'et\'e limite} (et nous notons $\, {\cal M} $) l'ensemble des
\'etats qui satisfont les conditions ``que l'on d\'esire obtenir'' \`a la fronti\`ere.
Ainsi, pour les cas habituels de fronti\`eres fluides, nous avons par exemple~: 

\smallskip \noindent  (4.3) $\qquad \displaystyle
 {\cal M} \,\,= \,\, \{ W_0 \} 	\qquad \qquad 	\qquad \qquad $ entr\'ee supersonique

\setbox11=\hbox{$ \,\,$ \'etats $\, W \,$ tels que les variables non conservatives }
\setbox12=\hbox {$ \,\, (\rho ,\, u ,\,p ) \,\, $ v\'erifient $\,\, {1\over2} \, u^2
\,+\, h(\rho,\,p) \,\,= \,\,H \,,\,\, {{p}\over{\rho^{\gamma}}} \,=\, S  \,  $  }
\setbox13= \vbox {\halign{#&# \cr \box11  \cr \box12     \cr}}
\setbox14= \hbox{ $\vcenter {\box13} $}
\setbox15=\hbox{\noindent   $\,\,  \displaystyle    {\cal M} \,\,= \,\, \left\{ \box14
\right\}  $} 
\setbox16=\hbox { \qquad   entr\'ee subsonique d'enthalpie totale et d'entropie
impos\'ees.  }
\setbox40= \vbox {\halign{#&# \cr \box15  \cr \box16    \cr}}
\setbox41= \hbox{ $\vcenter {\box40} $}
\setbox44=\hbox{\noindent  (4.4) $\qquad \displaystyle  \left\{
\box41 \right. $}  
\smallskip \noindent $ \box44 $

\smallskip \noindent  (4.5) $\qquad \displaystyle
 {\cal M} \,\,= \,\, \{ W \,/ \,\, p \,\,= \,\, {\overline p} \,\} \qquad \qquad  $
sortie subsonique de pression impos\'ee 
 
\smallskip \noindent  (4.6) $\qquad \displaystyle
 {\cal M} \,\,= \,\, \{ W \,/ \,\,u \,+\, c \,\leq \,0 \, \,\}
\qquad $ sortie supersonique.

\smallskip \noindent
La vari\'et\'e $\,  {\cal M} \,$ (\`a bord \'eventuellement,  {\it c.f.} (4.6)) est de
codimension $p.$ Notons que le cas d'une paroi solide d\'efinit \'egalement une
vari\'et\'e $\, {\cal M} ,\,$ qui d\'epend du pas de temps~: 

\smallskip \noindent  (4.7) $\qquad \displaystyle
 {\cal M} \,\,= \,\, \{ $ \'etat miroir de $\, W_1^n \,,\,\,$ d\'efini par la
relation (3.16) $\, \} \,. \, $

\bigskip \noindent   $\bullet \qquad  \,\,\, $  
A la fronti\`ere gauche ($x=0$), nous posons un {\bf ``probl\`eme de Riemann partiel''}
$\,\,P( {\cal M},\, W_1^n) \,\,$ entre la vari\'et\'e fronti\`ere $\, {\cal M} \,$
et l'\'etat dans la premi\`ere cellule $\, W_1^n .\,$ Cette notion, introduite dans
[Du87] et [DL89] ne correspond pas \`a un
probl\`eme de Cauchy comme pour le probl\`eme de Riemann usuel $\,\, R(W_g,\,W_d) . \,
\,$ Une solution  de $\,\,P( {\cal M},\, W_1^n) \,\,$ n'est tout d'abord d\'efinie que
dans l'espace des \'etats et consiste en une suite d'au plus codim$ \, {\cal M} \,$
$j$-ondes s\'epar\'ees par des \'etats constants, en suivant la d\'emarche classique
pour r\'esoudre la probl\`eme de Riemann (voir Lax [La73] par exemple.) On cherche
un \'etat $\, W = W_0 \, $ appartenant \`a $\,  {\cal M} ,\,$  et $p$ \'etats
interm\'ediaires $\, W ,\, \cdots \, W_{p-1} \,$ de sorte que 

\setbox15=\hbox {$ \,\, W \,=\, W_0 \in {\cal M} \,\, \,$ et $\,\,\,$  il existe 
$\,\, W_1 \,,\, W_2 \,, \cdots ,\, W_{p-1} \,\, (p\,=\, {\rm codim} {\cal M})\,\,$ }
\setbox16=\hbox {$ \,\,$tels que  $\,\,    W_1 \in U_{3-(p-1)}(W_0) \,  \,,\,\,  W_2
\in U_{3-(p-2)}(W_1) \,,\, \cdots \,,\, $ }
\setbox17=\hbox {$ \qquad  W_p \,=\, W_1^n \in U_3(W_{p-1})  \,. \,$  }
\setbox40= \vbox {\halign{#&# \cr \box15  \cr \box16     \cr \box17   \cr}}
\setbox41= \hbox{ $\vcenter {\box40} $}
\setbox44=\hbox{\noindent  (4.8) $\,\,\, \displaystyle  \left\{
\box41 \right. $}  
\smallskip \noindent $ \box44 $

\smallskip \noindent
Une fois l'\'etat $\, W \,$ d\'etermin\'e gr\^ace aux relations (4.8) (qui ne
consuisent pas toujours \`a un \'etat unique, voir par exemple le cas (4.4) dans
[Du87]), la solution de $\,\,P( {\cal M},\, W_1^n) \,\,$ est d\'efinie dans l'espace
$\, (x,\,t) \,$ comme la solution du probl\`eme de Riemann {\bf classique}
$\,\,\,R( W,\, W_1^n) \,.\,$ Le flux fronti\`ere $\, f_{1/2} \,$ est alors
simplement le flux num\'erique associ\'e au probl\`eme de Riemann 
$\,\,R( W,\, W_1^n) \,$: 

\smallskip \noindent  (4.9) $\qquad \displaystyle
f_{1/2} \,\,=\,\, \Phi \bigl( W ,\, W_{j\!+\!1}^n \bigr) \,, \qquad \qquad W \in
{\cal M} \,\,$ solution de (4.8). 

\smallskip \noindent
Comme les ondes du probl\`eme de Riemann peuvent prendre une c\'el\'erit\'e
arbitraire, on n'a pas forc\'ement $\,\, f_{1/2} \,=\, f(W) \,$ comme le proposent
implicitement Hemker et Spekreijse [HS86] par exemple. 

\bigskip \noindent   $\bullet \qquad  \,\,\, $  
Le cas d'une sortie subsonique de pression impos\'ee a \'et\'e abord\'e au cours de
la secon\-de partie. Les relations (3.9) sont dans ce cas \'equivalentes aux conditions
(4.8) et il suffit d'interpr\'eter les ondes du probl\`eme de Riemann $\,\,R(W,\,
W_1^n ) .\,\,$ Sur la figure 9 par exemple, l'\'etat $\, W_1^n \,$ correspond \`a une
sortie supersonique et la 3-d\'etente est enti\`erement localis\'ee dans le quadrant 
$\,\, \{ {{x}\over{t}}  \leq 0 \,,\, t \geq 0 \, \} .\,$  Le flux paroi est seulement
d\'etermin\'e par l'\'etat $\, W_1^n \,$: $\,\, f_{1/2} \,=\, f\bigl( W_1^n \bigr)
.\,$ [On comparera avec les autres cas de figure (3.10) et (3.12) relatifs \`a ce
m\^eme probl\`eme] .

%%%%%%%%%%%%%%%%%%%%%%%%%%%%%%%%      figure 2-9  %%%%%%%%%%%%%%%%%%%%%%%%%%%%%%%%    
%%%%   \bigskip   \smallskip 
%%%%   %%%  \centerline {  \epsfysize=4,0cm  \epsfbox  {../cnslim/Cnslim.fig9.epsf} } 
%%%%   \smallskip  \smallskip
% fin de la version linux
%%%%   \centerline { {\bf Figure 9}	\quad Sortie subsonique de pression impos\'ee, }
%%%%   \centerline { formul\'ee \`a l'aide du probl\`eme de Riemann partiel $\,\,P( {\cal
%%%%   M},\, W_1^n) \,\,$ }
%%%%   \centerline { et utilisant le flux mutivalu\'e de Osher. }    \smallskip 
%%%%%%%%%%%%%%%%%%%%%%%%%%%%%%%%%%%%%%%%%%%%%%%%%%%%%%%%%%%%%%%%%%%%%%%%%%%%%%%%%%%%% 

%%%%%%%%%%%%%%%%%%%%%%%%%%%%%%%%%%%%%%%%%%%%%%%%%%%%%%%%%%%%%%%%%%%%%%%%%%%%%%%%%%% figure  9
\bigskip
\centerline  {\includegraphics[width=.65\textwidth]   {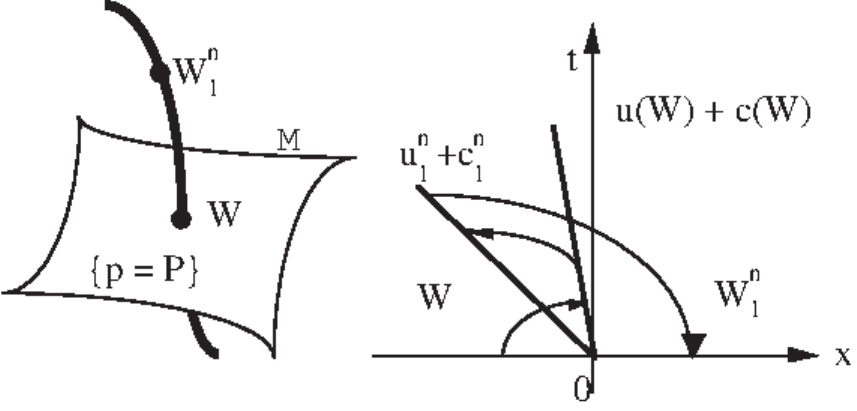}}
\smallskip  
\noindent {{\bf Figure 9}. Sortie subsonique de pression impos\'ee,
formul\'ee \`a l'aide du probl\`eme de Riemann partiel $\,\,P( {\cal M},\, W_1^n) \,\,$
et utilisant le flux mutivalu\'e de Osher.}
\bigskip
%%%%%%%%%%%%%%%%%%%%%%%%%%%%%%%%%%%%%%%%%%%%%%%%%%%%%%%%%%%%%%%%%%%%%%%%%%%%%%%%%%%

\bigskip \noindent   $\bullet \qquad  \,\,\, $  
Nous d\'etaillons maintenant les deux cas triviaux pour l'\'etude lin\'earis\'ee, \`a
savoir l'entr\'ee et la sortie supersonique. Dans le premier cas, la vari\'et\'e $\,
{\cal M} \,$ est un singleton ({\it c.f.} (4.3)) et le probl\`eme de Riemann partiel est en
fait un probl\`eme de Riemann classique. Ainsi, le flux paroi peut \^etre diff\'erent
de $\, f(W_g) ,\,$ comme les figures 1 \`a 3 le montrent. Dans le second cas, $\,
{\cal M} \,$ est une vari\'et\'e \`a bord ({\it c.f.}  (4.6)) et l'\'etat $\, W \,$ d\'efini
par la relation (4.8) est exactement $\, W_1^n\,$ si ce dernier appartient \`a $\,
{\cal M} ,\,$ mais est confondu avec l'\'etat sonique $\, W^* \,$ appartenant \`a la
3-d\'etente dans le cas contraire. Nous insistons sur le fait que le flux de sortie
{\bf ne se r\'eduit pas toujours \`a une extrapolation} du type 

\smallskip \noindent  (4.10) $\qquad \displaystyle
f_{1/2} \,\, = \,\, f\bigl( W_1^n \bigr) \qquad \qquad ( W_1^n  \in {\cal M}, $
d\'efinie en (4.6)) 

\smallskip \noindent 
comme nous l'avons propos\'e jusqu'ici, puisqu'il faut \'egalement tenir compte de
l'\'eventualit\'e 

\smallskip \noindent  (4.11) $\qquad \displaystyle
f_{1/2} \,\, = \,\, f\bigl( W^*  \bigr) \qquad \qquad ( W_1^n  \not\in {\cal M}, $
 {\it c.f.}  (4.6)) .

\bigskip \noindent   $\bullet \qquad  \,\,\, $  
La m\'ethode d\'ecrite dans ce paragraphe permet le calcul du flux paroi comme le
flux d'un probl\`eme de Riemann partiel, pos\'e entre la ``vari\'et\'e limite'' $\,
{\cal M}\,$ qui d\'ecrit les conditions physiques \`a imposer et l'\'etat $\, W_1^n
\,$ dans la cellule touchant le bord. Lorsqu'on utilise une r\'esolution {\bf exacte}
du probl\`eme de Riemann, nous g\'en\'eralisons la formulation de la condition limite
propos\'ee dans [DL87] et l'\'etat \`a la fronti\`ere peut \^etre
tr\`es {\bf \'eloign\'e} \`a la fois de l'\'etat $\, W_1^n \,$ et de la vari\'et\'e
$\, {\cal M} .\,$ Nous avons ainsi {\bf affaibli} la notion de condition limite pour
prendre en compte les  fortes ondes non lin\'eaires pr\'esentes \`a la fronti\`ere du
domaine de calcul. Les calculs alg\'ebriques relatifs aux cas de figure (4.3) \`a
(4.6) sont expos\'es dans [Du87] lorsqu'on utilise le solveur propos\'e par Osher
pour r\'esoudre le probl\`eme de Riemann. 

%% \newpage 
\bigskip \noindent {\bf 4.2)  \quad  	 Cas des faibles non lin\'earit\'es}

\noindent   $\bullet \qquad  \,\,\, $  
Nous supposons dans ce paragraphe que l'interaction au bord du domaine est faible,
{\it i.e.}

\smallskip \noindent  (4.12) $\qquad \displaystyle
W_1^n \,\, $ appartient \`a un voisinage de $\, {\cal M} \,.\, $

%%%%%%%%%%%%%%%%%%%%%%%%%%%%%%%%      figure 2-10    %%%%%%%%%%%%%%%%%%%%%%%%%%%%%%%%    
%%%%%  \bigskip    \smallskip 
%%%  \centerline {  \epsfysize=3,5cm  \epsfbox  {../cnslim/Cnslim.fig10.epsf} } 
%%%%%  \smallskip  \smallskip
% fin de la version linux
%%%%%  \centerline { {\bf Figure 10}	\quad  Entr\'ee subsonique faiblement non lin\'eaire. }
%%%%%   \smallskip 
%%%%%%%%%%%%%%%%%%%%%%%%%%%%%%%%%%%%%%%%%%%%%%%%%%%%%%%%%%%%%%%%%%%%%%%%%%%%%%%%%%%%%%%%% 

%%%%%%%%%%%%%%%%%%%%%%%%%%%%%%%%%%%%%%%%%%%%%%%%%%%%%%%%%%%%%%%%%%%%%%%%%%%%%%%%%%% figure  10
\bigskip
\centerline  {\includegraphics[width=.55\textwidth]   {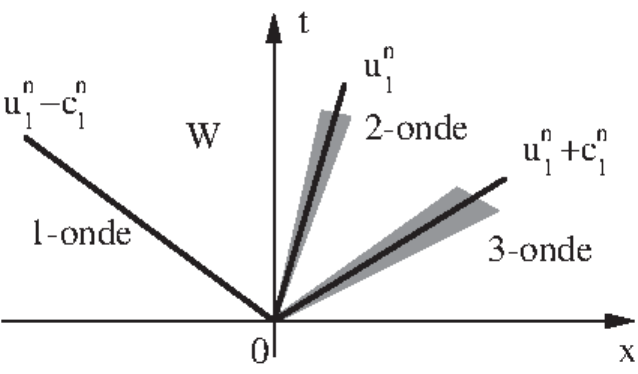}}
\smallskip  
\noindent {{\bf Figure 10}.  Entr\'ee subsonique faiblement non lin\'eaire.}
\bigskip
%%%%%%%%%%%%%%%%%%%%%%%%%%%%%%%%%%%%%%%%%%%%%%%%%%%%%%%%%%%%%%%%%%%%%%%%%%%%%%%%%%%

\smallskip  \smallskip \smallskip  \noindent
Quitte \`a prendre un voisinage assez petit, il est alors clair que le probl\`eme de
Riemann partiel a ses codim $\, {\cal M} \,$ ondes \`a l'int\'erieur du quadrant 
$\,\, \{ {{x}\over{t}}  > 0 \,,\, t \geq 0 \, \} \,$ (figure 10) puisque les vitesses
des ondes du  probl\`eme de Riemann partiel sont des grandeurs positives voisines des
valeurs propres $\, \lambda_j \bigl( W_1^n \bigr) .\,$ Le quadrant $\,\, \{
{{x}\over{t}}  \leq 0 \,,\, t \geq 0 \, \} \,$ contient donc seulement l'\'etat
constant $\, W \,$ d\'etermin\'e aux relations (4.8). On peut alors parler d'un
\'etat fronti\`ere et l'on a clairement 

\smallskip \noindent  (4.13) $\qquad \displaystyle
f_{1/2} \,\,= \,\, f(W) \,.\,$ 

\smallskip \noindent 
Le calcul de l'\'etat $\, W \,$ est effectu\'e soit par r\'esolution d'un probl\`eme
de Riemann (voir Osher-Chakravarthy [OC83] et Hemker-Spekreijse [HS86]), soit par
approximations des ondes non lin\'eaires par des discontinuit\'es de vitesses $\,
\lambda_j \bigl( W_1^n
\bigr) .\,$ On retrouve alors les relations de compatibilit\'e (relations (3.3)) qui,
jointes
\`a le condition 

\smallskip \noindent  (4.14) $\qquad \displaystyle
W \, \in \, {\cal M} \,$

\smallskip \noindent
permet de calculer l'\'etat fronti\`ere  de fa\c{c}on approch\'ee (voir
Veuillot et Viviand [VV78], Chakravarthy [Ch83], Cambier-Escande-Veuillot [CEV86]). 

%%%%%%%%%%%%%%%%%%%%%%%%%%%%%%%%%%%%%%%%%%%%%%%%%%%%%%%%%%%%%%%%%%%%%%%%%%%%%%%%%%%%%%%%%%%%%
\bigskip   
\noindent {\bf 4.3)  \quad  	 Etude d'un cas test}

%%%%%%%%%%%%%%%%%%%%%%%%%%%%%%%%      figure 2-11    %%%%%%%%%%%%%%%%%%%%%%%%%%%%%%%%     
%%%  \bigskip \centerline {  \epsfysize=10cm  \epsfbox  {../cnslim/Cnslim.fig11.epsf} } 
%%%%%%%  \smallskip  \smallskip
%%%%%%%  \centerline { {\bf Figure 11}	\quad 	Evolution de la vitesse en $\,  x=0 \, $ et 
%%%%%%%  $\,  x = 1 \, $  }
%%%%%%%  \centerline { pour le maillage de 80 cellules. }
%%%%%%%%%%%%%%%%%%%%%%%%%%%%%%%%%%%%%%%%%%%%%%%%%%%%%%%%%%%%%%%%%%%%%%%%%%%%%%%%%%%%%    

%%%%%%%%%%%%%%%%%%%%%%%%%%%%%%%%%%%%%%%%%%%%%%%%%%%%%%%%%%%%%%%%%%%%%%%%%%%%%%%%%%% figure  11
\bigskip
%% \centerline  {\includegraphics[width=.55\textwidth]   {fig-11.pdf}}
\centerline  {\includegraphics[width=.45\textwidth]   {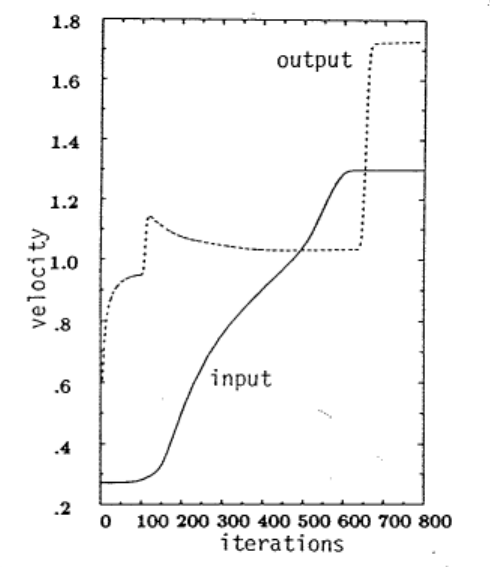}}

\smallskip  
\noindent {{\bf Figure 11}. Tuyère divergente. \'Evolution de la vitesse en $\,  x=0 \, $ (ligne continue)
et $\,  x = 1 \, $ (ligne pointill\'ee) pour le maillage contenant 80 cellules.}
\bigskip
%%%%%%%%%%%%%%%%%%%%%%%%%%%%%%%%%%%%%%%%%%%%%%%%%%%%%%%%%%%%%%%%%%%%%%%%%%%%%%%%%%%

\bigskip \noindent   $\bullet \qquad  \,\,\, $  
Nous avons \'etudi\'e un cas test tr\`es simple, d\'ej\`a abord\'e entre autres par
Yee, Beam et Warming [YBW82]~: il s'agit du calcul d'un \'ecoulement enti\`erement
supersonique dans une tuy\`ere divergente. La section $\, A(x) \,$ est donn\'ee par
la relation 

\smallskip \noindent  (4.15) $\qquad \displaystyle
A(x) \,\,= \,\, 1.598 \,+\, 0.347 \,\, {\rm th} \, (8x-4) \qquad 0 \, \leq \, x \, \leq
\, 1 \,$

\smallskip \noindent 
et l'entr\'ee supersonique est d\'efinie par l'\'etat suivant~: 

\smallskip \noindent  (4.16) $\qquad \displaystyle
\rho_g \,\,= \,\, 0.502 \,,\, \qquad u_g \, = \, 1.299 \,,\ \qquad p_g \,=\, 0.381
\,.$

\bigskip \noindent 
L'\'ecoulement supersonique recherch\'e est solution stationnaire du mod\`ele
quasi-monodimen\-sionnel des tuy\`eres, obtenu en adjoignant aux \'equations d'Euler
(2.1) un terme source afin de prendre en compte les variations de section dans les
bilans (voir par exemple Liu [Li82]). Nous avons effectu\'e le calcul [Du87] \`a
l'aide d'un sch\'ema explicite en temps (CFL=0.9), du premier ordre en espace et
divers maillages de 20, 40 et 80 cellules. Le probl\`eme de Riemann aux interfaces a
\'et\'e trait\'e \`a l'aide du flux d'Osher. Les conditions initiales correspondent
\`a un \'etat de vitesse nulle qui a m\^eme entropie et m\^eme enthalpie totale que
l'\'etat d\'ecrit en (4.16), {\it i.e.} 

\smallskip \noindent  (4.17) $\qquad \displaystyle
\rho_i \,\,= \,\,1 \,,\, \qquad u_i \, = \, 0 \,,\ \qquad p_i \,=\, 1 \,. \, $

\smallskip \noindent
Rappelons le traitement des conditions aux limites~: l'entr\'ee est supersonique donc
le flux limite est \'evalu\'e \`a l'aide du probl\`eme de Riemann entre l'\'etat
amont (4.16) et celui pr\'esent au pas de temps \'etudi\'e dans la premi\`ere
cellule. La sortie, supersonique \'egalement, est calcul\'ee \`a l'aide de l'un des
flux (4.10) ou (4.11) selon que l'\'etat de la derni\`ere cellule est supersonique
sortant ou non. 

%%%%%%%%%%%%%%%%%%%%%%%%%%%%%%%%      figure 2-12    %%%%%%%%%%%%%%%%%%%%%%%%%%%%%%%%      
%%%  \bigskip \centerline {  \epsfysize=10cm  \epsfbox  {../cnslim/Cnslim.fig12.epsf} } 
%%%%%%  \smallskip  \smallskip
%%%%%%  \centerline { {\bf Figure 12}	\quad 		Tuy\`ere divergente, 20 points de grille, 
%%%%%%  sortie supersonique.  }
%%%%%%  \centerline {  Evolution de la vitesse en $\,  x=0  \,$ }
%%%%%%  \centerline {  avec la condition limite traditionnelle   (4.18) (courbe 1)  }
%%%%%%  \centerline {  et avec la condition (4.9) fond\'ee sur le probl\`eme de Riemann
%%%%%%  (courbe 2). } 
%%%%%%%%%%%%%%%%%%%%%%%%%%%%%%%%%%%%%%%%%%%%%%%%%%%%%%%%%%%%%%%%%%%%%%%%%%%%%%%%%%%%%%  

%%%%%%%%%%%%%%%%%%%%%%%%%%%%%%%%%%%%%%%%%%%%%%%%%%%%%%%%%%%%%%%%%%%%%%%%%%%%%%%%%%% figure  12
\bigskip
%%% \centerline  {\includegraphics[width=.55\textwidth]   {fig-12.pdf}}
\centerline  {\includegraphics[width=.45\textwidth]  {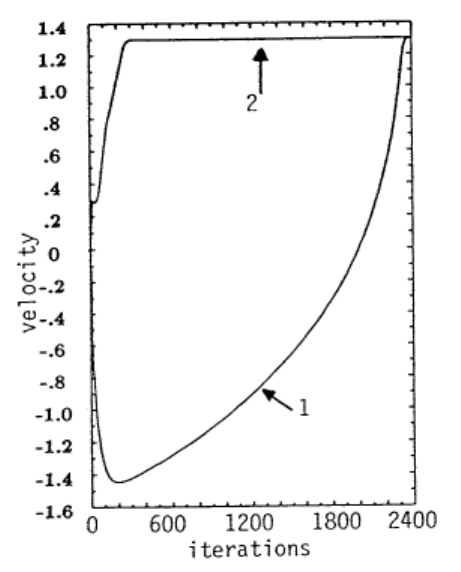}}
\smallskip  
\noindent {{\bf Figure 12}.
Tuy\`ere divergente, 20 points de grille, 
sortie supersonique. \'Evolution de la vitesse en $\,  x=0  \,$
avec la condition limite traditionnelle   (4.18) (courbe 1)
et avec la condition (4.9) fond\'ee sur le probl\`eme de Riemann (courbe 2).}
\bigskip
%%%%%%%%%%%%%%%%%%%%%%%%%%%%%%%%%%%%%%%%%%%%%%%%%%%%%%%%%%%%%%%%%%%%%%%%%%%%%%%%%%%

\bigskip \noindent   $\bullet \qquad  \,\,\, $  
Malgr\'e  l'inad\'equation entre une  condition initiale d'\'etat subsonique (!)
(4.17) et les conditions aux limites (4.16), le r\'egime stationnaire est ateint apr\`es une
\'evolution instationnaire importante mais r\'eguli\`ere (figure 11). Afin de mesurer
l'importance du choix de chacune des deux conditions au bord, nous avons effectu\'e
deux tests compl\'ementaires. Dans le premier cas, toutes choses \'egales par
ailleurs, nous changeons  le traitement num\'erique de la condition d'entr\'ee, en
rempla\c{c}ant le probl\`eme de Riemann $\, R(W_g \,,\, W_1^n) \,$ par une
\'evaluation ``classique'' du flux~: 

\smallskip \noindent  (4.18) $\qquad \displaystyle
f_{1/2} \,\,= \,\, f(W_g) \,.\,$

\smallskip \noindent
La convergence, qui demandait 200 pas de temps environ avec 20 points de grille, en
demande trois fois plus (figure 12) et la vitesse en $\, x=0 \,$ commence par \^etre
n\'egative (l'entr\'ee supersonique se comporte comme une sortie pendant la plus
longue partie du transitoire~!) avant d'atteindre la valeur finale exacte. Dans le
second cas, nous avons simplement rempla\c{c}\'e (par rapport \`a l'exp\'erience
initiale de la figure 11), la condition de sortie supersonique par une extrapolation
classique (4.10) \`a {\bf tous} les pas de temps. Le r\'esultat est encore plus
surprenant~: apr\`es une phase transitoire tr\`es br\`eve (trois fois plus rapide que
pour le test initial), on obtient une solution converg\'ee {\bf subsonique} (figure
13), et celle-ci est stable par raffinements successifs du maillage. La prise en
compte des fortes non-lin\'earit\'es de fronti\`ere est donc essentielle pour qu'une
\'evolution instationnaire converge vers ``la'' solution d\'efinie par les conditions
aux limites que l'on se donne, et ce avec une vitesse de convergence la plus
\'elev\'ee possible.

%%%%%%%%%%%%%%%%%%%%%%%%%%%%%%%%      figure 2-13    %%%%%%%%%%%%%%%%%%%%%%%%%%%%%%%%       
%%%  \bigskip \centerline {  \epsfysize=8,2cm  \epsfbox  {../cnslim/Cnslim.fig13.epsf} } 
%%%%%   \smallskip  \smallskip
%%%%%   \centerline { {\bf Figure 13}	\quad Tuy\`ere divergente, 80 points de grille.  }
%%%%%   \centerline {  Nombre de Mach \`a  convergence (courbe du bas) }
%%%%%   \centerline {  avec une condition limite de sortie calcul\'ee gr\^ace \`a
%%%%%   l'extrapolation (4.10) }
%%%%%   \centerline {  et solution exacte (courbe du haut). }  
%%%%%%%%%%%%%%%%%%%%%%%%%%%%%%%%%%%%%%%%%%%%%%%%%%%%%%%%%%%%%%%%%%%%%%%%%%%%%%%%%%%%%%%%% 

%%%%%%%%%%%%%%%%%%%%%%%%%%%%%%%%%%%%%%%%%%%%%%%%%%%%%%%%%%%%%%%%%%%%%%%%%%%%%%%%%%% figure  13
\smallskip %% \bigskip
%% \centerline  {\includegraphics[width=.50\textwidth]   {fig-13.pdf}}
\centerline  {\includegraphics[width=.65\textwidth]   {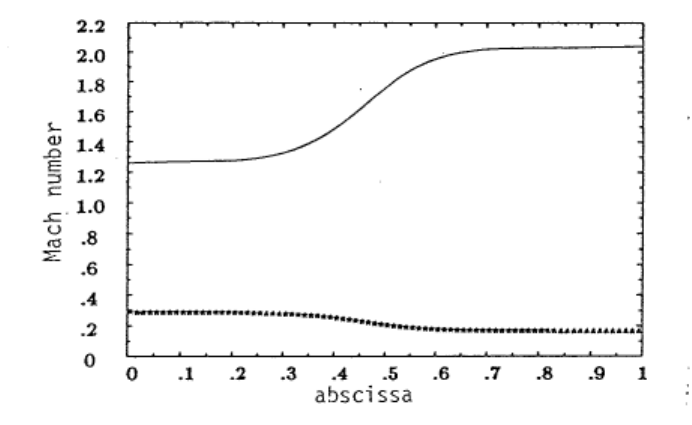}}
\smallskip  
\noindent {{\bf Figure 13}.
  Tuy\`ere divergente, 80 points de grille. Nombre de Mach \`a  convergence (\'etoiles)
  avec une condition limite de sortie calcul\'ee gr\^ace \`a l'extrapolation (4.10)
  et solution exacte (ligne continue).}
%% \bigskip
%%%%%%%%%%%%%%%%%%%%%%%%%%%%%%%%%%%%%%%%%%%%%%%%%%%%%%%%%%%%%%%%%%%%%%%%%%%%%%%%%%%

\bigskip \noindent   $\bullet \qquad  \,\,\, $  
La g\'eom\'etrie d'une tuy\`ere \'etant assez simple, et sans faire une analyse
d\'etaill\'ee de l'interaction des ondes instationnaires de choc et de d\'etente 
avec les ondes stationnaires induites par la g\'eom\'etrie (pour laquelle nous
renvoyons par exemple \`a Liu [Li82]), nous pouvons tenter d'interpr\'eter simplement
les divers ph\'enom\`enes qui se produisent au cours de ces trois exp\'eriences
num\'eriques.  Avec nos conditions limites initiales, nous impososns la r\'esolution
du probl\`eme de Riemann $\, R(W_g,\, W_1^n) ,\,$ lequel comporte un 1-choc de vitesse
n\'egative et un 3-choc de vitesse positive (figure 10). Par ailleurs, la condition
de sortie fait entrer une 1-d\'etente par l'aval du domaine de calcul. L'\'ecoulement
d'entr\'ee se stabilise vite autour d'une valeur interm\'ediaire relative au
probl\`eme de Riemann  $\, R(W_g,\, W_i) \, \,$ (figure 11) tandis que le 3-choc
interagit avec la 1-d\'etente au milieu de la tuy\`ere. Le r\'esultat de cette
interaction conduit \`a une augmentation r\'eguli\`ere de la vitesse en entr\'ee et
\`a un 1-choc qui sort de la tuy\`ere apr\`es 650 pas de temps environ (figure 11,
ligne fine).  %%%%%%%   pointill\'ee). 

\bigskip \noindent   $\bullet \qquad  \,\,\, $  
Lorsqu'on force la condition de flux (4.18), le 1-choc pr\'esent au temps $\, t=0 \,$
\`a l'entr\'ee est anim\'e d'une vitesse positive et est r\'efl\'echi brutalement par
la tuy\`ere (figure 12). La 1-d\'etente centr\'ee en $\, x=1 \,$ envahit alors
progressivement la tuy\`ere, mais les interactions internes sont plus fortes que dans
le cas pr\'ec\'edent et le ph\'enom\`ene met plus de temps \`a s'\'etablir. Si l'on
traite classiquement la condition de sortie, on ne cr\'ee pas de 1-d\'etente en
$\,x=1 \,$ et la tuy\`ere r\'eagit comme un tube \`a choc~: la 3-d\'etente sort du
domaine \'etudi\'e et la solution stationnaire est essentiellement li\'ee \`a
l'\'etat interm\'ediaire du probl\`eme de Riemann pos\'e \`a l'entr\'ee (figure 13). 

\bigskip \bigskip 
%%%%%%%%%%%%%%%%%%%%%%%%%%%%%%%%%%%%%%%%%%%%%%%%%%%%%%%%%%%%%%%%%%%%%%%%%%%%%%%  section 5
\noindent {\bf \large    5) \quad Conclusion} 
%%%%%%%%%%%%%%%%%%%%%%%%%%%%%%%%%%%%%%%%%%%%%%%%%%%%%%%%%%%%%%%%%%%%%%%%%%%%%%%%%%%%%%%%%%

\noindent   $\bullet \qquad  \,\,\, $  
Nous avons pass\'e en revue les diverses approches math\'ematiques et num\'eri\-ques
pour l'\'etude du probl\`eme des conditions aux limites associ\'ees \`a la
r\'esolution des \'equations d'Euler de la dynamique des gaz. La prise en compte de
fortes non-lin\'earit\'es peut s'\'ecrire th\'eoriquement \`a l'aide d'une 
{\bf in\'egalit\'e d'entro\-pie \`a la fronti\`ere}. Lorsqu'on interpr\`ete cette derni\`ere
gr\^ace au {\bf pro\-bl\`eme de Riemann}, les effets non lin\'eaires pr\'esents \`a la
fronti\`ere sont pris en compte naturellement avec la {\bf m\'ethode des volumes
finis} qui {\bf affaiblit} la condition limite~; la discussion classique sur les
diff\'erents cas de figure (entr\'ee ou sortie, sub ou supersonique) se
r\'einterpr\`ete \`a l'aide d'un {\bf probl\`eme de Riemann partiel} entre une
vari\'et\'e limite et un \'etat fluide. Les premiers tests num\'eriques montrent
l'int\'er\^et d'uen telle m\'ethode.

\bigskip \noindent   $\bullet \qquad  \,\,\, $  
Toutefois, il reste de nombreuses questions sans r\'eponse actuellement~: quelle est
la stabilit\'e non lin\'eaire des sch\'emas num\'eriques ainsi d\'efinis ? Comment
aborder le cas d'un domaine multidimensionnel ? Enfin, dans la plupart des
applications, on \'etudie les fortes non-lin\'earit\'es {\bf \`a l'int\'erieur} du
domaine de calcul et la fronti\`ere ne joue pas un r\^ole d\'eterminant pour la prise
en compte des effets non lin\'eaires~!

%% \newpage
%%%%%%%%%%%%%%%%%%%%%%%%%%%%%%%%%%%%%%%%%%%%%%%%%%%%%%%%%%%%%%%%%%%%%%%%%%%%  references
\bigskip \bigskip  \bigskip      \noindent {\bf  \large  References }
%%%%%%%%%%%%%%%%%%%%%%%%%%%%%%%%%%%%%%%%%%%%%%%%%%%%%%%%%%%%%%%%%%%%%%%%%%%%%%%%%%%%%%%%%%

\smallskip

\smallskip \hangindent=9mm \hangafter=1 \noindent
[Au84]  $\,$ J. Audounet, ``Solutions discontinues
param\'etriques des syst\`emes de lois de conservation et des probl\`emes aux limites
associ\'es'', {\it S\'eminaire}, Universit\'e Toulouse 3, 1983-84.

\smallskip \hangindent=9mm \hangafter=1 \noindent
[BLN79] C. Bardos, A.Y. Leroux, J.C. N\'ed\'elec, ``First Order Quasilinear Equations with Boundary Conditions'',
{\it Communications in Partial Differential Equations}, volume~4, pages~1017-1034, 1979.

\smallskip \hangindent=9mm \hangafter=1 \noindent
[BS87]  $\,$  A. Benabdallah, D. Serre, ``Probl\`emes aux 
limites pour les syst\`emes hyperboliques non-lin\'eaires de deux \'equations \`a une
dimension d'espace'', {\it Comptes Rendus de l'Aca\-d\'emie des Sciences}, Paris, tome~303,
S\'erie~1, pages~677-680, 1987.

\smallskip \hangindent=9mm \hangafter=1 \noindent
[BS77]  $\,$ J.S. Bramley, D.M. Sloan, ``A comparison of
boundary methods for the numerical solution of hyperbolic systems of equations'',
{\it Journal of Engineering Mathematics},  volume~11, n$^{\rm o}$3, pages~227-239, 1977.

\smallskip \hangindent=9mm \hangafter=1 \noindent
[CEV86] $\,$  L. Cambier, B. Escande, J.P. Veuillot, ``Calcul d'\'ecoulements internes
\`a grand nombre de Reynolds par r\'esolution des
\'equations de Navier-Stokes,   {\it La Recherche A\'erospatiale},  n$^{\rm o}$1986-6,
pages~415-432, 1986.

\smallskip \hangindent=9mm \hangafter=1 \noindent
[Ch83]  $\,$ S. Chakravarthy, ``Euler solutions, Implicit 
schemes and Boundary Conditions'', {\it AIAA Journal}, volume~21,  n$^{\rm o}$5, pages~699-706, 1983. 

\smallskip \hangindent=9mm \hangafter=1 \noindent
[CG85]  $\,$ P. Colella, H.M. Glaz, ``Efficient Solution
Algorithms for the Riemann Problem for Real Gases'', {\it Journal of Computational Physics},
volume~59, pages~264-289, 1985.

\smallskip \hangindent=9mm \hangafter=1 \noindent
[CF48] $\,$  R. Courant, K.O. Friedrichs. {\it 
Supersonic Flow and Shock Waves,}  Interscience, New-York, 1948.

\smallskip \hangindent=9mm \hangafter=1 \noindent
[DP83] $\,$  R. Di Perna, ``Convergence of the Viscosity
Method for Isentropic Gas Dynamics'', {\it Communications on Pure and Applied Mathematics},
volume~91, pages~1-30, 1983. 

\smallskip \hangindent=9mm \hangafter=1 \noindent
[Du87] $\,$ F. Dubois, ``Boundary Conditions and the
Osher Scheme for the Euler Equations of Gas Dynamics'', {\it Rapport interne
n$^{\rm o}$170 du Centre de Math\'ematiques Appliqu\'ees de l'\'Ecole Polytechnique}, Palaiseau,
septembre 1987.

\smallskip \hangindent=9mm \hangafter=1 \noindent
 [DL87]  $\,$ F. Dubois, P. Le Floch,  ``Condition \`a
la limite pour un syst\`eme de lois de conservation'',  {\it Comptes Rendus de
l'Acad\'emie  des Sciences}, Paris, tome~304, S\'erie~1, pages~75-78, 1987.

\smallskip \hangindent=9mm \hangafter=1 \noindent
 [DL88]  $\,$  F. Dubois, P. Le Floch, ``Boundary
Conditions for Nonlinear Hyperbolic Systems of Conservation Laws'', {\it Journal of
Differential Equations}, volume~71, n$^{\rm o}$1, pages~93-122, 1988.

\smallskip \hangindent=9mm \hangafter=1 \noindent
 [DL89]  $\,$  F. Dubois, P. Le Floch, ``Boundary 
Conditions for Nonlinear Hyperbolic Systems of Conservation laws'', Second International
Conference on Hyperbolic Problems (J. Ballmann, R. Jeltsch Editors), %%  Josef Ballmann, Rolf Jeltsch
{\it Notes on Numerical Fluid Dynamics}, volume~24, pages~96-104, Vieweg, Braunschweig, 1989.

\smallskip \hangindent=9mm \hangafter=1 \noindent
  [DG88]  $\,$  B. Dubroca, G. Gallice, ``Probl\`eme mixte
pour un syst\`eme de lois de conservation monodimensionnel'',  {\it Comptes Rendus de
l'Acad\'emie  des Sciences},  Paris, tome~306, S\'erie~1, pages~317-320, 1988. 

\smallskip \hangindent=9mm \hangafter=1 \noindent
[GB53] $\,$  P. Germain, R. Bader, ``Unicit\'e des
\'ecoulements avec chocs dans la m\'ecanique de Burgers'', {\it Note technique ONERA
OA}  n$^{\rm o}$11/1711-1, mai 1953. 

\smallskip \hangindent=9mm \hangafter=1 \noindent
[Gl65] $\,$ J. Glimm,  ``Solutions in the Large for
Nonlinear Hyperbolic Systems of Conservation Laws'',
{\it Communications on Pure and Applied Mathematics}, volume~18, pages~95-105, 1965.

\smallskip \hangindent=9mm \hangafter=1 \noindent
[Go59]  $\,$ S.K. Godunov, ``A Difference Method for the
Numerical Computation of Discontinuous Solutions of the Equations of Fluid Dynamics'',
{\it Matematicheskii Sbornik}, volume~47, pages~271-290, 1959.

\smallskip \hangindent=9mm \hangafter=1 \noindent
[Go61] $\,$ S.K. Godunov, ``An intersting class of
quasilinear Systems'', {\it  Doklady Akademii Nauk SSSR}, volume~139, pages~521-523~;
voir aussi  {\it Soviet mathematics - doklady}, American Mathematical Society, volume~2, pages~947-949, 1961.

\smallskip \hangindent=9mm \hangafter=1 \noindent
[GZIKP79]   S.K. Godunov, A. Zabrodine, M. Ivanov,
A. Kraiko, G. Prokopov. {\it  R\'esolution num\'erique des probl\`emes
multidimensionnels de la dynamique des gaz},  \'Editions de Moscou, 1979.

\smallskip \hangindent=9mm \hangafter=1 \noindent
[GT87] $\,$ M. Goldberg, E. Tadmor, ``Convenient
Stability Criteria for Difference Approximations of Hyperbolic Initial-Boundary Value
Problems II'', {\it Mathematics of Computation}, volume~48,   n$^{\rm o}$178, pages~503-520, 1987.

\smallskip \hangindent=9mm \hangafter=1 \noindent
[Gu85] $\,$ B. Gustafsson, ``Numerical Boundary
Conditions, in Large Ccale Computations in Fluid Mechanics'' (Engquist, Osher,
Somerville Editors), {\it Lectures in Applied Mathematics}, volume~22, AMS, Providence,
pages~279-308, 1985.

\smallskip \hangindent=9mm \hangafter=1 \noindent
[HLV83] A. Harten, P.D. Lax, V. Van Leer,   
``On Upstream Differencing and Godunov-type Schemes for Hyperbolic Conservation Laws'',
{\it  SIAM Review,} volume~25, n$^{\rm o}$1, pages~35-61, janvier 1983.

\smallskip \hangindent=9mm \hangafter=1 \noindent
[HS86] $\,$  P.W. Hemker, S.P. Spekreijse, ``Mutiple Grid
and Osher's Scheme for the Efficient Solution of the Steady Euler Equations'',
{\it Applied Numerical  Mathematics}, volume~2, pages~475-493, 1986.

\smallskip \hangindent=9mm \hangafter=1 \noindent
[Hi86] $\,$  R.L. Higdon. Initial-Boundary Value Problems
for Linear Hyperbolic Systems, {\it SIAM Review}, volume~28,  n$^{\rm o}$2, pages~177-217,
1986. 

\smallskip \hangindent=9mm \hangafter=1 \noindent
[Ke71] $\,$ C.P. Kentzer, ``Discretization of Boundary
Conditions on Moving Discontinuities'', in {\it Lecture Notes in Physics} (M. Holt Ed.),
volume~8, Springer Verlag, Berlin, pages~108-113, 1971.

\smallskip \hangindent=9mm \hangafter=1 \noindent
[Kr70] $\,$  H.O. Kreiss. Initial Boundary Value Problems
for Hyperbolic Systems, {\it Communications on  Pure and Applied Mathematics}, volume~23, pages~277-298, 1970.

\smallskip \hangindent=9mm \hangafter=1 \noindent
[Kv70] $\,$  S. Kru$\breve {\rm z} $kov, ``First Order  Quasi-Linear Systems
in Several Independant Variables'', {\it Matematicheskii Sbornik}, volume~123, pages~228-255,
et {\it Mathematics of the USSR-Sbornik}, volume~10,  n$^{\rm o}$2, pages~217-243, 1970. 

\smallskip \hangindent=9mm \hangafter=1 \noindent
[LL54] $\,$  L. Landau, E. Lifchitz.  {\it Fluid Mechanics},
1954, Editions de Moscou, 1967.

\smallskip \hangindent=9mm \hangafter=1 \noindent
[La71] $\,$  P.D. Lax, ``Shock Waves and Entropy'', in
{\it Contribubutions to Nonlinear Functional Analysis} (Zarantonello Ed.),
Academic Press, New York, pages~603-634, 1971. 

\smallskip \hangindent=9mm \hangafter=1 \noindent
[La73] $\,$  P.D. Lax.  {\it Hyperbolic Systems of Conservation
Laws and the Mathematical Theory of Shock Waves},   Conference Board in Mathematical
Sciences, volume~11, Society for Industrial and Applied Mathematics, Philadelphia, 1973.

\smallskip \hangindent=9mm \hangafter=1 \noindent
[LW60] $\,$   P.D. Lax, B. Wendroff, ``Systems of
Conservation Laws'',  {\it Communications on  Pure and Applied Mathematics}, volume~13, pages~217-237, 1960.

\smallskip \hangindent=9mm \hangafter=1 \noindent
[LF88]  $\,$ P. Le Floch, ``Explicit Formula for Scalar
Nonlinear Conservation Laws with Boundary Conditions'',
{\it Mathematical  Methods in the Applied Sciences}, volume~10, pages~265-287, 1988.

\smallskip \hangindent=9mm \hangafter=1 \noindent
[LN88] $\,$ P. Le Floch, J.C. N\'ed\'elec, ``Explicit Formula for Weighted Scalar
Nonlinear Hyperbolic Conservation Laws'', {\it Transactions of the American Mathematical Society},
volu\-me~308, pages~667-683, 1988.

\smallskip \hangindent=9mm \hangafter=1 \noindent
[Le81]  $\,$ A. Lerat. {\it Sur le calcul des solutions
faibles des syst\`emes hyperboliques de lois de conservation \`a l'aide de sch\'emas
aux diff\'erences}, Th\`ese d'Etat, Universit\'e Paris 6, 1981.

\smallskip \hangindent=9mm \hangafter=1 \noindent 
[Li77]    $\,$ T.P. Liu, ``Initial-Boundary Value Problems
for Gas Dynamics'', {\it Archive for Rational Mechanics and Analysis}, volume~64,
pages~137-168, 1977.

\smallskip \hangindent=9mm \hangafter=1 \noindent 
[Li82]    $\,$  T.P. Liu, ``Transonic Gas Flow in a Duct of
Vatying Area'',   {\it Archive for Rational Mechanics and Analysis}, volume~80,  n$^{\rm o}$1,
pages~1-18, 1982. 

\smallskip \hangindent=9mm \hangafter=1 \noindent 
[Mc69]    $\,$  R.W. Mac Cormack, ``The Effect of Viscosity
in Hypervelocity Impact Cratering'', {\it AIAA Paper}  n$^{\rm o}$69-354, 1969.

\smallskip \hangindent=9mm \hangafter=1 \noindent 
[MO75]   $\,$    A. Majda, S. Osher, ``Initial-Boundary Value
Problems for Hyperbolic Equations with Uniformly Characteristic Boundaries'', 
{\it Communications on  Pure and Applied Mathematics}, volume~28, pages~607-675, 1975.

\smallskip \hangindent=9mm \hangafter=1 \noindent 
[MBGB87] P. Mazet, F. Bourdel, R. Greborio, J. Bor\'ee, ``Application de la m\'ethode
variationnelle d'entropie \`a la r\'esolution des \'equations d'Euler'',
{\it ONERA, Centre d'Etudes et de Recherches de Toulouse}, Rapport Interne, 1987.

\smallskip \hangindent=9mm \hangafter=1 \noindent 
[Mo68]  $\,$  G. Moretti, ``The Importance of Boundary
Conditions in the Numerical Treatment of Hyperbolic Equations'', {\it Polytechnic
Institute of Brooklyn}, PIBAL Report n$^{\rm o}$68-34, 1968 ; voir aussi
{\it Physics of  Fluids}, volume~12, pages~II-13–II-20, 1969.

\smallskip \hangindent=9mm \hangafter=1 \noindent 
[NS77]    $\,$  T. Nishida, J. Smoller, ``Mixed Problems for
Nonlinear Conservation Laws'', {\it Journal of Differential Equations},
volume~23, pages~244-269, 1977. 

\smallskip \hangindent=9mm \hangafter=1 \noindent 
[Ol57]   $\,$  O. Oleinik, ``Discontinuous Solutions of
Nonlinear Differential Equations'', {\it Uspehi Matemati$\breve {c} $eskih Nauk}, volume~12, pages~3-73~; %% ${\check{\it c}}$  {\v{c}}
voir aussi {\it American Mathematical Society Translations: Series 2}, volume~26, pages~95-172, 1957.

\smallskip \hangindent=9mm \hangafter=1 \noindent 
[OS78]   $\,$ J. Oliger, A. Sundstr\"om, ``Theoretical and Practical
Aspects of Some Initial Boundary Value Problems in Fluid Dynamics'', 
{\it SIAM Journal on Applied Mathematics}, volu\-me~35, pages~419-446, 1978. 

\smallskip \hangindent=9mm \hangafter=1 \noindent 
[Os81]   $\,$ S. Osher, ``Solution of Singular Perturbation
Problems and Hyperbolic Systems of Conservation Laws'', in 
{\it Mathematical Studies} n$^{\rm o}$47 (Axelsson-Franck-Van der Sluis Editors), pages~179-205,
North Holland, Amsterdam, 1981.

\smallskip \hangindent=9mm \hangafter=1 \noindent 
[Os84]   $\,$ S. Osher, ``Riemann Solvers,
the Entropy Condition and Difference Approximations'', {\it SIAM Journal of Numerical
Analysis},  volume~21, pages~217-235, 1984.

\smallskip \hangindent=9mm \hangafter=1 \noindent 
[OC83]   $\,$  S. Osher, S. Chakravarthy, ``Upwind Schemes and
Boundary Conditions with Applications to Euler Equations in General Geometries'',
{\it Journal of Computational Physics}, volume~50, pages~447-481, 1983.

\smallskip \hangindent=9mm \hangafter=1 \noindent 
[RM67]    $\,$  R.D. Richtmyer, K.W. Morton. {\it Difference
Methods for Initial-Value Problems}, Interscience Publishing, J. Wiley \& Sons, New
York, 1967.

\smallskip \hangindent=9mm \hangafter=1 \noindent 
 [Ri81]  $\,$ A. Rizzi, ``Computation of Rotational
Transonic Flow'', in {\it Numerical Methods for the Computation of Inviscid Transonic Flows
with Shock Waves}, (Rizzi-Viviand Editors), Vieweg Verlag, Braunschweig, pages~153-161, 1981.

\smallskip \hangindent=9mm \hangafter=1 \noindent 
[Ro72]  $\,$  P.J. Roache. {\it Computational Fluid
Dynamics}, Hermosa Publishers, Albukerque, 1972.

\smallskip \hangindent=9mm \hangafter=1 \noindent 
[Ro81]  $\,$ P. Roe, ``Approximate Riemann Solvers,
Parameter Vectors and Difference Schemes'', {\it Journal of Computational Physics},
volume~43, pages~357-372, 1981.

\smallskip \hangindent=9mm \hangafter=1 \noindent 
[Sm83]    $\,$  J. Smoller.  {\it Shock Waves and
Reaction-Diffusion Equations}, Springer Verlag, Berlin, 1983.

\smallskip \hangindent=9mm \hangafter=1 \noindent 
[St84]    $\,$ B. Stoufflet. {\it R\'esolution num\'erique des
\'equations d'Euler des fluides parfaits compressibles par des sch\'emas implicites
en \'el\'ements finis}, Th\`ese de Docteur-Ing\'enieur, Universit\'e Paris 6,
1984.

\smallskip \hangindent=9mm \hangafter=1 \noindent 
 [Th87]  $\,$ K.W. Thomson, ``Time-Dependent Boundary
 Conditions for Hyperbolic Systems'',  {\it Journal of Computational Physics}, volume~68,
 pages~1-24, 1987.

\smallskip \hangindent=9mm \hangafter=1 \noindent 
 [VL79]   $\,$   B. Van Leer, ``Towards the Ultimate
Conservative Difference Scheme V. A Second Order Sequel to Godunov's Method'',
{\it Journal of Computational Physics}, volume~32, n$^{\rm o}$1, pages~101-136, 1979. 

\smallskip \hangindent=9mm \hangafter=1 \noindent 
[VL84]   $\,$   B. Van Leer, ``On the relation between the
Upwind-Differencing Schemes of Godunov, Engquist-Osher and Roe'',
{\it SIAM Journal on Scientific Computing}, volume~5,  n$^{\rm o}$1, pages~1-20, 1984.

\smallskip \hangindent=9mm \hangafter=1 \noindent 
[VV78]  $\,$  H. Viviand, J.P. Veuillot, ``M\'ethodes
pseudo-instationnaires pour le  calcul d'\'ecoule\-ments transsoniques'', {\it Publication
ONERA}  n$^{\rm o}$1978-4, 1978.

\smallskip \hangindent=9mm \hangafter=1 \noindent 
[YBW82]    $\,$ H. Yee, R. Beam, R. Warming, ``Boundary
Approximations for Implicit Schemes for One-Dimensional Inviscid Equations of Gas
Dynamics'', {\it AIAA Journal}, volume~20, pages~1203-1211, 1982.

\end{document}